\documentclass[11pt]{amsart}
\usepackage[english]{babel}
\usepackage{graphics}
\usepackage{epsfig}
\usepackage{latexsym,amssymb,amsfonts,amsmath, amscd, euscript}

\usepackage{enumerate}

\newcommand{\N}{\mathbb{N}}

\newcommand{\Q}{\mathbb{Q}}

\newcommand{\B}{\mathcal B}

\newtheorem{definition}{{\bf Definition}}

\newtheorem{theorem}{{\bf Theorem}}
\newtheorem{proposition}{\noindent {\bf Proposition}}
\newtheorem{corollary}{\noindent {\bf Corollary}}

\newtheorem{claim}{\noindent {\bf Claim}}

\newtheorem{question}{\noindent {\bf Question}}
\newtheorem{questions}[question]{\noindent {\bf Questions}}
\newtheorem{example}{\noindent {\bf Example}}

\newtheorem{problem}{\noindent{\bf Problem}}

\newtheorem{remark}{\noindent {\bf Remark}}

\def\proofref #1 {{\noindent  {\bf Proof} (#1).}\ }

\newtheorem{lemma}[definition]{\noindent {\bf Lemma}}
\def\endproof{\hfill {\kern 6pt\penalty 500
\raise -0pt\hbox{\vrule \vbox to5pt {\hrule width 5pt
\vfill\hrule}\vrule}}}
\def\centerpicture #1 by #2 (#3){\leavevmode
        \vbox to #2{
        \hrule width #1 height 0pt depth 0pt
        \vfill
        \special{pictfile #3}}}
\begin{document}
\title{On  scattered posets with finite dimension}\dedicatory{
To the memory of Eric C.Milner (1928-1997)}

\author{Maurice Pouzet}\address{Math\'ematiques, ICJ, Universit\'e Claude-Bernard Lyon1,
43 Bd 11 Novembre 1918, 69622 Villeurbanne Cedex, France}
\email{pouzet@univ-lyon1.fr}
\author{Hamza Si Kaddour}
\address{Math\'ematiques, ICJ, Universit\'e Claude-Bernard Lyon1,
43 Bd 11 Novembre 1918, 69622 Villeurbanne Cedex, France}
\email{sikad@univ-lyon1.fr}
\author{Nejib Zaguia}
\address{SITE, Universit\'e d'Ottawa, 800 King Edward Ave, Ottawa, Ontario, K1N6N5,Canada}
\email{zaguia@site.uottawa.ca}

\thanks{Research done under the auspices of the CMCU Franco-Tunisien  "Outils math\'ematiques pour l'informatique". This research was completed while the authors visited each others. Support provided by the university of Ottawa is gratefully acknowledged}



\date{\today }

\begin{abstract}         We discuss a possible  characterization, by means of forbidden configurations, of posets which are embeddable in a product of
finitely many scattered chains.  \end{abstract}

 \maketitle

 \section*{Introduction and presentation of the results}
 
 A fundamental result, due to Szpilrajn \cite{szpilrajn30},  states that every order on a set  is the intersection of  a family of linear orders on this set. The \emph{dimension} of the order, also called the dimension of the ordered set,   is then defined as the minimum cardinality of such a family (Dushnik, Miller \cite{dushnik-miller}). Specialization of Szpilrajn's result to several types of orders have been studied \cite{{bonn-pouz2}}. An  ordered set (in short poset), or its order, is \emph {scattered} if it does not contain a  subset which is ordered as the chain $\eta$ of rational numbers. Bonnet and Pouzet \cite{bonn-pouz1} proved that   \emph{a poset is scattered if and  only if the order is the intersection of scattered linear orders}. 
 It turns out that there are  scattered posets whose order is the  intersection of finitely many linear orders but which cannot be the intersection of finitely many scattered linear orders. We give nine examples in Theorem \ref{thm:main}. This naturally leads to the following question:
\begin{question} \label{que:scatt} If  an order  is the intersection of finitely many linear scattered orders, does this order the intersection of $n$ many scattered linear orders, where $n$ is the dimension of this order?
 \end{question}
 
 We do not have the answer even for dimension two orders.  We cannot even answer this: \begin{question} If an order of dimension two  is the intersection of three scattered linear orders, does this order the intersection of two scattered linear orders?
 \end{question}
 Question \ref{que:scatt} is a special instance of the following general qestion:
 
 \emph{Given a positive integer $n$, which orders are intersection of  at most $n$ scattered linear orders?}

We propose an approach based on the notion of obstruction. 

Let $n$ be a non negative  integer; let  $\mathcal {L}(n)$, resp. $\mathcal {L}_{\mathcal S}(n)$   be the class of posets $P$ whose order is the the  intersection of at most $n$ linear orders, resp. at  most $n$ scattered linear orders.   Set $\mathcal {L}(<\omega):=\bigcup_{n<\omega}\mathcal {L}(n)$ and  $
\mathcal {L}_{\mathcal S}(<\omega):=\bigcup_{n<\omega}\mathcal {L}_{\mathcal S}(n)$.

These four classes are \emph{closed under embeddability}, that is if $\mathcal C$ is one of these classes, then for every poset $P\in \mathcal C$, a  poset $Q$ belongs to $\mathcal C$ whenever it  is embeddable in  $P$ (that is $Q$ is isomorphic to an induced subposet of $P$). Say that an \emph {obstruction} to a class $\mathcal C$ as above is any poset not belonging to $\mathcal C$. Then  such a class $\mathcal C$ can be characterized by  obstructions, eg  as the class of posets in which no obstruction to $\mathcal C$ is embeddable. But, it can be also characterized by means of smaller collections of obstructions. If $\mathcal B$ is a class of poset,  denote by  $Forb(\mathcal B)$  the class of posets  in which no member of $\mathcal \B$ is embeddable. 

With this terminology, we may ask:

\emph{Find  $\mathcal B$ as simple as possible such that $\mathcal {L}_{\mathcal S}(n)=Forb(\mathcal B)$.}

The following question emerges immediately:
\begin{question} Is there  a cardinal $\lambda$ such that every obstruction to $\mathcal {L}_{\mathcal S}(n)$ contains an obstruction  of size at most $\lambda$? \end{question}

As it can be easily seen, the existence of such a cardinal for an arbitrary  class closed under embeddability follows readily from the \emph{Vo$\check{p}$enka principle}, a  strong set theoretical principle which could be inconsistent with usual  set theoretical axioms. It implies  the existence of large cardinal numbers (eg supercompact cardinals) and its consistency is  implied by the existence of huge cardinals  (see \cite{jechst78} pp. 413--415). 

In the case  of $\mathcal {L}_{\mathcal S}(n)$ we do not know if $\lambda$ exists. In fact, \emph{we conjecture that it exists and is countable.}

The same general question for $\mathcal {L}(n)$ has a simpler answer: each obstruction contains a finite one. Indeed, as it is well known,  \emph{ a poset $P$ belongs to $\mathcal {L}(n)$ whenever for every \emph {finite} subset  $A$ of $P$ the poset induced by $P$ on $A$ is also in $\mathcal {L}(n)$} (this striking fact is a consequence of the compactness theorem of first order logic - for a proof, see the survey \cite{kelly}).  Furthermore, if $\mathcal Crit(\mathcal {L}(n))$ denotes the collection of minimal obstructions (that is the collection of finite posets $Q$ whose dimension is larger than $n$, whereas every proper subposet has dimension at most $n$), then $\mathcal {L}(n)=Forb(\mathcal Crit(\mathcal {L}(n)))$.  Members of $\mathcal Crit(\mathcal {L}(n))$ have dimension  $n+1$;  these posets are the so-called   \emph{$n+1$-irreducible posets} \cite{trotter}. For $n=1$, there is just one: the two element antichain. For $n=2$, a  complete description has been given by D.Kelly in 1972 (see \cite{kelly}). For $n>2$ a description seems to be hopeless;  in fact, the problem to decide whether or not a finite poset belongs to $\mathcal {L}(n)$ is NP-complete. If $\mathcal C=\mathcal {L}(<\omega)$,  every obstruction contains  a countable one (this easily follows from the finitary result mentionned above), hence $\mathcal {L}(<\omega)=Forb(\mathcal B)$ where $\mathcal B$ is a set of countable posets, each with a countable dimension.  
In terms of obstructions, Question \ref{que:scatt} amounts to:
\begin{question} Is $Crit(\mathcal {L}(n))$ determines  $\mathcal {L}_{\mathcal S}(n)$  within $\mathcal {L}_{\mathcal S}(<\omega)$?\end{question}

\begin{figure}[htbp]
\centering
\includegraphics[width=3in]{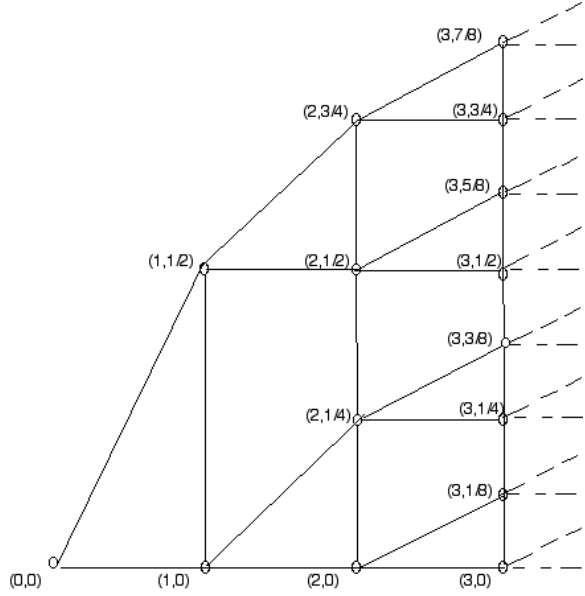}
\caption{$ \Omega(\eta)$}
\label{Omega}
\end{figure}

We rather consider the following:
\begin{question}
Is $\mathcal {L}_{\mathcal S}(<\omega)$ can be determined within  $\mathcal {L}(<\omega)$ by a finite set $\mathcal B_{\mathcal S}$ of obstructions?
  \end{question}
  We provide ten examples of  obstructions. All are countable and have dimension at most $3$. 
 In  order to present these examples, 
 we denote by $P^*$ the dual of  a poset $P$, we denote by $\check P$ the set $P$ equipped with the strict order $<$. We denote by 
  $B(\check P)$ the poset defined as
follows: the underlying set is $P\times \{0,1\}$,  the ordering defined
by $(x,i)  < (y,j)$ if  $i < j$ and $x< y$. This poset is the \emph{open split} of $P$.  It is clearly bipartite, moreover $B(\check P^*)$ is order-isomorphic to $B(\check P)^*$. Let $T_2$  be the infinite binary tree and let $\Omega (\eta )$ be the infinite binary tree in which each level is
totally ordered by an increasing way from the left to the right (see Figure \ref{Omega} for an equivalent representation).

We prove: 

 \begin{theorem} \label {thm:main}  A poset whose the order is the intersection of finitely many scattered linear orders contains no isomorphic copy of $\eta$, $T_2$,  $\Omega (\eta)$, $B(\check \eta)$, $B(\check T_2)$,  $B(\check \Omega (\eta))$ and their dual. 
\end{theorem}
 
 Since $\eta$ and $B(\check \eta)$ are self dual, this list contains only ten members. In fact, these members do not embed in each other (Lemma \ref{lem:incomp}).

 \begin{problem} Is this  list determines the class $\mathcal {L}_{\mathcal S}(<\omega)$ of orders  which are intersection of finitely many scattered linear orders within the class $\mathcal {L}(<\omega)$ of orders which are intersection of finitely many linear orders? 
\end{problem}

 The reader will notice that each of our obstructions distinct from $\eta$ contains an infinite antichain. This is general. Indeed, if a poset  $P$ is scattered with no infinite antichain,  each linear extension of the order on $P$ is scattered (\cite {bonn-pouz1}, see also \cite {bonn-pouz2}), hence if $dim (P)=n$, the order is the intersection of $n$ scattered linear orders. 
 
 The occurence of open splits in Theorem  \ref{thm:main} asks for an explanation.  We present one, despite the fact that it is not fully satisfactory. It is based on the notion of split  rather than open split. If $P$ is a poset, the \emph{split} of $P$ is the poset  $B(P)$  whose underlying set is $P\times\{0, 1\}$ ordered by:
$$(x,i)<(y,j)\;  \text{if } \; x\leq y\; \text{ and}  \; i<j.$$

We prove:
\begin{theorem}\label{bip2} Let $P$ be  a poset. Then $P\in  \mathcal L_{\mathcal S}(<\omega)$ if and only if $B(P)\in \mathcal L_{\mathcal S}(<\omega)$.
\end{theorem}

The analogous equivalence with $B(\check P)$ instead of $B(P)$ is in general false. But, if  $P$ is $\eta$, $T_2$,  $\Omega (\eta)$ or their dual, $B(P)$ and $B(\check P)$ can be embedded in each other (Lemma \ref{lem:comp}).  Hence, in order to prove  Theorem \ref{thm:main} it suffices to prove that  $\eta$, $T_2$,  $\Omega (\eta)$ and their dual are obstructions to $\mathcal {L}_S(<\omega)$ and to apply Theorem \ref{bip2}. In order to do that, we introduce a peculiar object:
the topological closure $\overline {N(P)}$ in the powerset $\mathfrak P(P)$ of the MacNeille completion $N(P)$ of a poset $P$. As a poset, $\overline {N(P)}$ is an algebraic lattice. 

We prove:
 \begin{theorem}\label{thm:second}
 Let $P$ be a poset and $n$ be a positive integer.  Then the following properties are equivalent:
 \begin{enumerate}[{(i)}]
 \item The order on $P$ is the intersection of $n$ scattered linear orders;
 \item $\overline{N(P)}$ is embeddable  into a product of $n$ scattered linear orders.
 \end{enumerate}
 Moreover, if one of these conditions hold, $\overline {N(P)}$ is topologically scattered. 
 \end{theorem}
 With this result at hand, in order to show that if $P$ is $\eta$, $T_2$ or  $\Omega (\eta)$, $P$ is an obstruction, it suffices to observe that  $\overline{N (P )}$ is not topologically scattered. We give the proof of this fact in Section \ref{section:proofs}. 

Note that while $N(P)$ and $N(P^*)$ are dually isomorphic, $\overline{N(P)}$ and $\overline {N(P^*)}$ are not. Hence, one can be topologically scattered, whereas the other is not. For an example, $\overline{N(T_{2}^*)}$ is topologically scattered and $\overline{N(T_{2})}$ is not.  
\begin{question} If $dim (P)\leq n$ and both $\overline{N(P)}$ and $\overline {N(P^*)}$ are topologically scattered does the order on  $P$ is the intersection of $n$ scattered linear orders?\end{question}

After such  unsuccessful attempt of a description of $\mathcal L_{\mathcal S}(<\omega)$ by means of obstructions, we looked at  subclasses $\mathcal C$ of $\mathcal L_{\mathcal S}(<\omega)$ such that  every member of  $\mathcal L_{\mathcal S}(<\omega)$ can be embedded in a member of $\mathcal C$. It turns out that the class of
scattered distributive  lattices of finite dimension has this property. In fact: 
\begin{theorem}\label{thm:distributivelattice}Let $T$ be a distributive lattice. The following properties are equivalent:
\begin{enumerate}[{(i)}]
\item The order on $T$ is the intersection of $n$ scattered linear orders.

\item $T$ is isomorphic to a sublattice of a  product of $n$ scattered chains.

\item $dim(T)\leq n$ and $T$ is order-scattered.
 \end{enumerate}

\end{theorem}

We  also consider extensions of  our initial question.  

Instead of linear orders, we consider interval orders and instead of scattered linear orders,  interval orders which can be represented as intervals of a scattered chain. Instead of ordered sets we consider incidence structures,  we replace linear orders  by Ferrers relations, we replace MacNeille completion by Galois lattices  and scattered linear orders  by  Ferrers relations whose Galois lattice is scattered. We obtain an extension of Theorem \ref{thm:second} (see Theorem \ref {thm:extension}). From our study, it  follows that a positive answer to our initial question implies a positive answer to the extensions we consider.  The basic objets of our study are incidence structures and Galois lattices. One of our key result is a  property of the topological closure of Galois lattices (Theorem \ref
{thm:topologicalbouchet}) which refines Bouchet's Coding theorem (\cite {bouchetetat}, see also \cite {bouchet}, see Theorem \ref{bouchet}).

To conclude, we mention a   specialization of our question which comes  
from the following observation.  All finite ordered sets of dimension $2$ are obtained as follows.
 Let $\underline n:= \{0,2,...,n-1\}, n \geq 2$ and  $C$
be the linear ordering $0<1<\dots <n-1$ on $\underline n$. Let $\sigma$ be a
permutation on $\underline n$, distinct of the identity map. Define the
order $\leq _{\sigma}$
on $\underline n$ by $x\leq _{\sigma} y$ if and only if $x \leq_C y$ and
$\sigma (x) \leq _C \sigma (y)$. Let $P_{\sigma}:= (\underline n ,\leq _{\sigma})$, and
$C_{\sigma}:= \{(x,y): 
\sigma (x) \leq _C \sigma (y)\}$, then $\leq _{\sigma}$ is the intersection of $C$ and $C_{\sigma}$.  Thus $P_{\sigma}$ has dimension $2$. For infinite posets, 
even countable, the situation is quite different.  Order which are 
intersection of two orders of type $\omega$ are close to finite orders.  A characterization in terms of obstructions is included in the
following:

\begin{theorem}
 An order on an infinite set  is the intersection of $n$ linear orders of type $\omega$
if and only if:
\begin{enumerate}[{(i)}]
\item  The order has dimension at most $n$.
\item The poset does not  contain an infinite antichain, an
infinite decreasing chain, the chain $\omega +1$ and the direct sum $\omega\oplus 1$ of the chain $\omega$ with the one-element chain.
\end{enumerate}
\end{theorem}

See Proposition 4.1 and Corollary 4.2 of \cite{pouzet-sauer}. 

 A more general question is the following:
 \begin{question} \label{que:sametype}
 Given a positive integer $n$ and an order type $\alpha$. Which orders are the intersection of $n$ linear orders of the same type $\alpha$? 
  \end{question}
 More specifically
 \begin{question}
 Characterize by means of obstructions the posets which are embeddable into posets whose order is the intersection of $n$ linear orders of the same type $\alpha$? 
 \end{question}

 This paper is composed as follows. Section  \ref{section:ingredients} contains the definitions of the main notions with a  development on incidence structures and Galois lattices and coding. It  includes our refinement  of  Bouchet's Coding theorem, and also some basic facts on Ferrers relations, interval orders and dimension. Section \ref {section:scattered} contains a discussion on the notions of scattered dimension, including  Theorem \ref {thm:extension}. Sections \ref{section:proofs} and \ref {section:proof4} contains the proofs of the results presented above. Section \ref{section:two-dimensional} contains a characterization  of orders which are intersection of two scattered linear orders (Theorem \ref{equiv}).

 \section{Ingredients} \label{section:ingredients}

 Our terminology follows \cite{daveypriestley} and \cite{gratzer98}.  Among  set theoretical notations, we point out that if  $f$ is a map from a set $E$ to a set $F$, and $A$ is a subset of $E$,   the set $\{f(x): x\in A\}$, the \emph{image} of $A$ by $f$, 
is denoted by  $f[A]$ rather than $f(A)$. \subsection{Order, lattices and topology}
As usual, a {\it poset} is the pair $P$ formed of a set $E$ and an order $\varepsilon$ on $E$.   If the order is \emph {linear}(or total), the poset is a \emph{chain}. The {\it dual} of $P:= (E,
\varepsilon)$ is
$P^*:=(E, \varepsilon ^{-1})$. If this causes no confusion, we will denote an order on $E$  by the symbol
$\leq$ and its complement by
$\not
\leq$;  we will denote the equality relation by $=$ (and, when needed, by $\Delta _E:=\{ (x,x): x\in E\}$),  we identify $P$ with $E$, writing $x\in P$ instead of $x\in E$. We will denote by $x \parallel _P y$ the fact that two elements $x$ and $y$ of $P$ are incomparable. Given a poset $P:=(E,
\leq)$, a subset $I$ of $E$ is an {\it initial 
segment} (or is {\it closed downward}) if $x\leq 
y$ and $y\in I$ imply $x\in I$. Let $X$ be a subset of $E$, we set:
\begin{align}\label{downarrow}
\downarrow\hskip -2pt  X:= \{y\in E: y\leq x \text{  for some  }  x\in X\}.
\end{align}
This set is an initial segment, in fact the least initial segment containing $X$.  We  say that $\downarrow\hskip-2pt X$ is generated by $X$. If $X$ contains only one element $x$,  we write $\downarrow\hskip -2pt  x$ instead of $\downarrow\hskip -2pt  \{x\}$. An initial segment of this form is \emph{principal}. We set $down(P):=\{\downarrow x: x\in P\}$.

We denote by ${\bf  I}(P)$ the set of initial segments of $P$ ordered by inclusion. For example, ${\bf I}((E, \Delta_E))={\mathfrak P}(E)$ the power set of
$E$ ordered by inclusion, whereas ${\bf I}((\Q,\leq))$ is the {\it Cantor  chain}. We also denote by ${\bf  I}_{<\omega}(P)$ the set of  finitely generated initial segments of $P$ ordered by inclusion. 
 An  {\it ideal} of  $P$  is a non empty initial segment  $I$ which is up-directed, that is every pair $x, y\in I$ has an upper bound $z\in I$.   We denote by ${\mathcal J}(P)$ the set of  ideals of $P$ and by $\mathcal J^{\neg \downarrow\hskip -2pt }(P)$ the subset  of non-principal ideals 
 of $P$. 
   Let $N(P)$ be the set made of intersections of principal  initial segments of $P$. Ordered by inclusion, $N(P)$ is a complete lattice, called the \emph{MacNeille completion of $P$}. 
    
A \emph{join-semilattice} is a poset $P$ such that every two elements $x,y$ have a least upper-bound, or join, denoted by $x\vee y$. If $P$ has a least element, that we denote $0$, this amounts to say that every finite subset of $P$ has a join.  An element $a$ in a lattice
$L$ is {\it compact}\index{compact} if for every $A\subset L$,
$a\leq \bigvee A$ implies $a\leq \bigvee A'$ for some finite
subset $A'$ of $A$. The lattice $L$ is {\it compactly
generated}\index{compactly generated} if every element is a
supremum of compact elements. A lattice is {\it
algebraic}\index{algebraic} if it is complete and compactly
generated. Algebraic lattices and join-semilattices with a least element  are  sides of the same coin.  Indeed, the set $K(L)$ of compact elements of an   algebraic lattice $L$ is a join-semilattice with a least element and $L$  is isomorphic to  the  set $\mathcal J(K(L))$  of ideals of $K(L)$, ordered by inclusion. Conversely, the  set $\mathcal J(P)$ of ideals of a join-semilattice $P$ having a least element, once ordered by inclusion, is an algebraic lattice, and the subset $K(\mathcal J(P))$ of its compact elements is isomorphic to $P$. 
We note that if $P$ is an arbitrary poset, ${\bf  I}(P)$ is an algebraic lattice and $K({\bf  I}(P))={\bf I}_{ <\omega}(P)$. Hence,  $\mathcal J({\bf I}_{ <\omega}(P))$ is order isomorphic to ${\bf I}(P)$. 
We also note that $\mathcal J(P)$  is
the set of join-irreducible elements of $I(P)$; moreover,
${\bf I}_{<\omega}(\mathcal J(P))$ is order-isomorphic to ${\bf I}(P)$ whenever $P$ has no infinite
antichain. 
  
Identifying  the power set 
 $\mathfrak{P}(E)$ of a set $E$ with $2^{E}$, we may view it  as a topological space.  A basis of open sets consists of  subsets  of  the form 
$O(F,G):=\{X\in \mathfrak{P}(E): F\subseteq X $ and $ G\cap X=\emptyset \}$, where $F, G$ are 
finite subsets of $E$.  As it is customary, we denote by $\overline {\mathcal F}$ the topological closure of a subset $\mathcal F$ of $\mathfrak P(E)$. 
 Recall that a  compact  totally disconnected space is called a {\it Stone space}, whereas a {\it Priestley space} is a set $X$ together with a  topology and an ordering which  is compact and \emph {totally order disconnected} in the sense that for every $x,y\in X$ such that $x\not \leq y$ there is some clopen initial segment containing $y$ and not $x$. Closed subspaces of $\mathfrak{P}(E)$, with the inclusion order  added, are Priestley spaces \cite{prie}. For an example, we recall that if  $L$ is an algebraic lattice then, with the topology induced by the product topology on $\mathcal J(K(L))$, it becomes a Priestley space. Priestley spaces are associated to bounded distributive lattices as Stone spaces are associated to Boolean algebras. We will recall in Section \ref{section:proof4} the properties we need about the relationship between Priestley spaces and distributive lattices. We refer to \cite{prie} and to \cite {daveypriestley} for an introduction to Stone-Priestley duality and to \cite{Gierz and all} for more on topologically ordered structures.
\subsection{Basic facts}
We will need the following basic result due to O.Ore and T.Hiraguchi (see \cite{schroder}): 
\begin{lemma}\label{lem:dimensionproduct}
Let $P$ be a poset and $\kappa$ be a cardinal. The order on $P$ is the intersection of $\kappa$ linear orders if and only if $P$ is embeddable in a product of $\kappa$ chains. 
\end{lemma}

\begin{lemma} \label{lem:product}
\begin{enumerate}
\item Let $P$ and $Q$ be two posets.
If $P$  is embeddable in  $Q$ then  $\mathcal J(P)$ is embeddable  in
$\mathcal J (Q)$.
\item Let $(P_i)_{i\in I}$ be a family of posets, then 
$\mathcal J(\Pi_{i\in I} P_i)$ is order-isomorphic to $\Pi_{i\in I} \mathcal J(P_i)$ provided that $I$ is finite.
\end{enumerate}
\end{lemma}
\begin{proof}
\noindent The proof  of Item 1 is immediate.  For Item 2,  let $A$ be a subset of $Q:=\Pi_{i\in I} P_i$. Given $i\in I$,  let $p_i: Q\rightarrow P_i$ be the $i$-th projection and  $p_i[A]$ be the image of $A$. Finally, set $\overline p(A):= (p_i[A])_{i\in I}$. We prove  that 
$\overline p$ induces an order-isomorphism from $\mathcal J(Q)$ onto $\Pi_{i\in I} \mathcal J(P_i)$. From this, Item 2 follows. Let $A\in \mathcal J(Q)$. First, we claim that  $\overline p(A)\in \Pi_{i\in I} \mathcal J(P_i)$. Indeed, let $i\in I$. Since $p_i$ is order-preserving and $A$ is up-directed, $p_i[A]$ is up-directed. Furthermore, $p_i[A]\in \mathbf {I}( P_i)$. Indeed,  let $x\in P_i$ and  $y\in p_i[A]$ such that $x\leq y$. Let $\overline y\in A$ such that $p_i(\overline y)=y$. Let $\overline x\in  \Pi_{i\in I} P_i$ defined by $\overline x_i:=x$ and $\overline x_j:=y_j$ for $j\not= i$. Then  $\overline x\leq \overline y$. Since $A\in \mathbf{I}(\Pi_{i\in I} P_i)$, $\overline x\in A$, and thus $x\in p_i[A]$, proving that $p_i[A]\in \mathbf{I}(P_i)$. Since $p_i[A]$ is up-directed, $p_i[A]\in \mathcal J(P_i)$. Thus  $\overline p(A)\in \Pi_{i\in I} \mathcal J(P_i)$ as claimed. Next, let $\overline A:=(A_i)_{i\in I}\in \Pi_{i\in I} \mathcal J(P_i)$. Then, trivially, $\Pi_{i\in I}A_i\in\mathcal J(Q)$.  Since all $A_i$'s are non-empty, $\overline p(\Pi_{i\in I}A_i)=\overline A\in\mathcal J(Q)$, proving that $\overline p$ is surjective. To conclude that $\overline p$ is an isomorphism, we note that  $A=\Pi_{i\in I}p_i[A]$ for every  $A\in \mathcal J(Q)$. Indeed, we have trivially $A\subseteq \Pi_{i\in I}p_i[A]$. For the reverse inclusion, let $\overline x\in \Pi_{i\in I}p_i[A]$. For each $i\in I$, select $\overline y(i)\in A$ such that  $\overline y(i)_i=\overline x_i$. Since $I$ is finite and $A$ is up-directed, there is $\overline z\in A$ which majorizes each  $\overline y(i)$. Due to our choices, $\overline z$ majorizes $\overline x$. Thus  $\overline x\in A$, as required.
\end{proof}

Let $E$ be  a set and  $\mathcal F$ be a subset of $\mathfrak P(E)$. We say that $\mathcal F$ is \emph{closed under intersections} if $\cap \mathcal F'\in \mathcal F$ for every 
 subset $\mathcal F'$ of $\mathcal F$ (with the convention that $\mathcal F'=E$ if $\mathcal F'=\emptyset$). We denote by $\mathcal F^{\wedge}$ the set of intersections of members of $\mathcal F$, (in particular $E\in \mathcal F^{\wedge}$. Hence $\mathcal F$ is closed under intersections if and only if $\mathcal F=\mathcal F^{\wedge}$.  Sets closed under intersections are usually called \emph{Moore families}.  As it is well known, a Moore family $\mathcal F$ is topologically closed in $\mathfrak P(E)$ if and only if it is closed under  unions of up-directed subfamilies. Moore families correspond 
to closure systems, those which are topology closed  to \emph{algebraic closure systems}  \cite{gratzer98}, \cite {Gierz and all}.  We will need the following fact.

 \begin{proposition}\label{prop:topoclo}

 Let $E$ be  a set and  $\mathcal F$ be a subset of $\mathfrak P(E)$. 
 Then $\overline{\mathcal F}^{\wedge}=\overline{\mathcal F^{\wedge}}$.
 \end{proposition} \begin{proof} It relies on the following claims.
 \begin{claim}\label {claim:wedge1}
 $\overline{\mathcal F^{\wedge}}$ is closed under intersections.
 \end{claim} 
 \noindent{\bf Proof of Claim \ref{claim:wedge1}.} Set $\mathcal G:=\overline{\mathcal F^{\wedge}}$. Let $\mathcal G'\subseteq \mathcal G$ and $X:= \bigcap \mathcal G'$. We prove that $X\in \mathcal G$. For that we prove that $O(F,G)\cap {\mathcal F^{\wedge}}\not = \emptyset$ for each finite $F\subseteq X $ and finite $G\subseteq E\setminus X$. We may suppose $X\not =E$ (otherwise, since   $E\in {\mathcal F^{\wedge}}$,  $X\in \mathcal G$ as required).  Let $ a\in G$.  Since $X= \bigcap \mathcal G'$ there is some $X_a\in \mathcal G'$ such that $X\subseteq X_a \subseteq E\setminus \{a\}$. Since $X_{a}\in \mathcal G$, there is some $Y_{F, a}\in O(F,\{a\})\cap \mathcal  F^{\wedge}$. Let $X_{F}:= \cap_{a\in G} X_{a}$. Clearly $X_F\in O(F,G)\cap\mathcal  F^{\wedge}$. This proves our claim
 \endproof
 \begin{claim}\label {claim:wedge2}
 $\overline{\mathcal F}^{\wedge}$ is topologically closed. \end{claim} 
 \noindent{\bf Proof of Claim \ref{claim:wedge2}.} Set $\mathcal G:=\overline{\mathcal F}^{\wedge}$. Let $X\in \overline {\mathcal G}$. Then $O(F,G)\cap \mathcal G\not = \emptyset$ for each finite $F\subseteq X $ and finite $ G\subseteq E\setminus X$. This implies that for every finite subset $F\subseteq X$, $a\in E\setminus X$, $O(F,\{a\})\cap \overline{\mathcal F}\not = \emptyset$. Let $ a\in E\setminus X$.  Since $\overline{\mathcal F}$ is compact, the intersection  $\bigcap \{ O(F,\{a\}): F\subseteq X, \;  F\;  \text{finite}\} \cap \mathcal \overline{\mathcal F}$ is non empty. Pick $X_a$ in this intersection. Let $X':=\bigcap_{a\in  E\setminus X} X_{a}$. Then $X'\in \overline{\mathcal F}^{\wedge}= \mathcal G$. But, since each $X_a$ contains $X$, $X=X'$, hence $X\in \mathcal G$. It follows that $\overline {\mathcal G}= \mathcal G$. This proves our claim. \endproof

 From Claim \ref{claim:wedge1} we deduce that $\overline{\mathcal F}^{\wedge}$ is included into $\overline{\mathcal F^{\wedge}}$ and from Claim \ref{claim:wedge2} the reverse inclusion. 
 \end{proof}
 
 If $P$ is a poset, we have  $N(P)= {down(P)}^{\wedge}$.  Hence,  Proposition \ref{prop:topoclo} yields immediately:
 \begin{corollary}
 $\overline {N(P)}=\overline {down(P)}^{\wedge}$.
 \end{corollary}

We recall the following fact (\cite {bekk-pouz-zhan} Corollary 2.4).
\begin{lemma}\label{lem:ideals1}
\begin{equation}\label{eq:ideals1}
down(P)\subseteq \mathcal {J}(P)\subseteq \overline {down (P)}\setminus\{\emptyset\}.
\end{equation}
 In
particular,   the topological closures in $\mathfrak{P}(P)$ of
$down(P)$  and $\mathcal {J}(P)$ are the same. 
\end{lemma}

\begin{lemma} \label{lem:ideals2}Let $P$ be  a join-semilattice with a least element. Then:
\begin{equation}
\overline {down(P)}=\overline {N(P)}= \mathcal J(P).
\end{equation}
\end{lemma}
\begin{proof} We start with the following:
\begin{claim}\label{claim:ideal2}\begin{equation}\label {eq:ideal2}
down(P)\subseteq  N(P)\subseteq \mathcal J(P).
\end{equation}
\end{claim}
\noindent{\bf Proof of Claim \ref{claim:ideal2}.}
Trivially, $down (P)\subseteq \mathcal J(P)$. Since $P$ is  a join-semilattice with a least element, $\mathcal J(P)$ is closed under intersection.  Hence, $N(P)$ which is made of the
 intersections  of members of $down(P)$ is included into $\mathcal J(P)$. \endproof

With Lemma \ref{lem:ideals1} this yields:\begin{equation}\label {eq:ideal3}
\overline {down(P)}\subseteq \overline {N(P)}\subseteq 
\overline{\mathcal J(P)}\subseteq \overline {down (P)}.
\end{equation}
To conclude, we note that $\mathcal J(P)$ is topologically closed. Indeed, $\mathcal J(P)$ is closed under union of up-directed sets and  as, observed  above,  it is closed under intersection. \end{proof}
\subsection{Incidence structures and coding} 
Let $E,F$ be two sets. A {\it binary relation from 
$E$ to $F$} is any subset $\rho$ of the cartesian 
product $E\times F$. As usual, we denote by 
$x\rho y$ the fact that 
$(x,y)\in \rho$ and by $x \neg \rho y$ the negation. The triple
$R:=(E,\rho ,F)$ is an {\it incidence structure}; its
 {\it complement} is $\neg R
:=(E, \neg \rho,F)$, where $\neg \rho:= E\times F\setminus \rho$, whereas its  {\it dual }
is $R^{-1} :=(F,\rho ^{-1},E)$, where $\rho ^{-1} := \{ (y,x) : (x,y) \in
\rho \}$. We set $L_R(Y):= \{ x\in E
: \{ x\} \rho Y \}$, resp. $U_R(X):= \{ y\in F
: X\rho \{ y\}  \}$, for each $Y\subseteq F$, resp. $X\subseteq E$. And we use  $L_R(y)$ and $U_R(x)$ for $L_R(\{y\})$ and $U_R(\{x\})$. With these notations, we have $U_R(X)=L_{R^{-1}}(X)$. The
sets $Gal(R):=\{ L_R(Y): Y\subseteq F\}$ and 
$Gal(R^{-1}):=\{ U_R(X): X\subseteq E\}$ are closed  under intersection; hence, once ordered by inclusion, they are 
complete lattices. Ordered by inclusion, $Gal(R)$ is the 
 {\it Galois lattice} of $R$.  A fundamental result is that $Gal(R^{-1})$ is isomorphic to $Gal(R)^*$, the {\it dual} of $Gal(R)$. If $P:= (E, \leq)$ is a poset,  $Gal((E, \not\geq ,E))= {\bf I}(P)$, whereas $Gal((E, \leq ,E))=N(P)$.
 
 Let $R:=(E,\rho ,F)$, $R':=(E',\rho ',F')$ be two incidence
structures, a {\it coding from} $R$ {\it into} $R'$ is a 
pair of maps  
$f: E\rightarrow E' , \ \ g: F\rightarrow F'$ such that 
 $$x\rho y \Longleftrightarrow f(x)\rho ' g(y)$$
  for all $x\in E$ and $y\in F$.  
When such a pair exists, we say that $R$ has  a {\it  coding into} $R'$.

\begin{example} If $R:= (E, \rho, F)$ is an incidence structure, the pair 
$(f,g)$,  where  $f(x):= L_R\circ U_R(x)$ for $x\in E$ and   $g(y):= L_R(y)$ for $y\in F$, is a coding from $R$ into $(Gal(R), \subseteq, Gal(R))$. \end{example}
 
 If $E=F$ and $E'=F'$, the pairs $(E,\rho)$, $(E',\rho')$ are
{\it binary relational structures} (or simply, {\it directed graphs}) and a map $f:E\rightarrow E'$ is an {\it embedding} if it is
one-to-one and 
$$x\rho y   \Longleftrightarrow f(x)\rho' f(y)$$ for all $x, y\in E$. When such a  map exists, we say that $(E,\rho)$ {\it is embeddable into} $(E', \rho')$. 

\begin{example}If  $\rho$ and $\rho'$ are two orders and $(E', \rho')$ is a complete lattice, $R$ has a coding into $R'$ if and only if $(E, \rho)$ is embeddable in $(E, \rho')$.
 \end{example}

Bouchet's Coding theorem (\cite {bouchetetat}, see also \cite {bouchet}) is a striking illustration of the links between coding and embedding. 

\begin {theorem}\label {bouchet}   Let $T$ be a  complete lattice and 
$R$ be an incidence structure, then $R$ has a coding into $(T,\leq  ,T)$ if
and only if $Gal(R)$ is embeddable in $T$. \end {theorem}   

\begin{corollary}\label{lem:coding-embedding}
Let  $R:= (E, \rho, F)$ and $R':= (E', \rho', F')$ be two incidence structures. Then  $Gal(R)$ is embeddable in $Gal(R')$ whenever $R$ has a coding into $R'$. 
\end{corollary}

We will need the following strengthening of Corollary \ref{lem:coding-embedding}.
\begin{theorem}\label{thm:topologicalbouchet}
Let  $R:= (E, \rho, F)$ and $R':= (E', \rho', F')$ be two incidence structures. If $R$ has a coding into $R'$ then there is an embedding $\phi$ from $\overline {Gal(R)}$ into $\overline {Gal(R')}$ and a continuous and order preserving map $\psi$ from a closed subspace $\mathcal H$ of $\overline {Gal(R')}$ onto $\overline{Gal(R)}$ such that $\psi\circ \phi =1_{\overline {Gal(R)}}$.
\end{theorem}
\begin{proof}
Let $(f,g)$ be a coding from $R$ into $R'$. Let $f^{d}: \mathfrak P(E')\rightarrow \mathfrak P(E)$ be defined by $f^d(X'):= f^{-1}(X')$ for $X'\subseteq E'$. The map $f^d$ is continuous. With the fact that $\mathfrak P(E')$ is compact, it follows that:
\begin{equation}\label{eq:close1}
f^{d}[\overline {\mathcal F'}]= \overline {f^d[\mathcal F']}
\end{equation}
for every $\mathcal F'\subseteq \mathfrak P(E')$.

Let $\mathcal F:= \{L_{R'}(g[Y]): Y\subseteq F\}$.
Clearly, $\mathcal F$ is closed under intersection and included into  $Gal(R')$. Furthermore, since $(f,g)$ is a coding:

\begin{equation}\label{eq:close2}
f^{d}(L_{R'}(g(y)))=L_R(y)
\end{equation}
for every $y\in F$.

Hence, \begin{equation}\label{eq:close3}
f^{d}(L_{R'} (g[Y]))=L_R(Y)
\end{equation}
for every $Y\subseteq  F$.

This implies: 
\begin{equation}\label{eq:close4}
f^{d}[ {\mathcal F}]=  Gal(R).
\end{equation}
With equation (\ref{eq:close1}), this yields:
\begin{equation}\label{eq:close5}
f^{d}[ \overline {\mathcal F}]=  \overline{Gal(R)}.
\end{equation}
Set $\mathcal H:= \overline {\mathcal F}$ and $\psi:= f^{d}_{\restriction \mathcal H}$. Clearly $\mathcal H$ is a closed subset of $Gal(R')$ and $\psi$ is a continuous and order preserving map from  $\mathcal H$ into $\overline {Gal(R)}$. Let $X\in \overline{ Gal(R)}$. According to equation (\ref{eq:close5}), $X$ 
belongs to the range of $\psi$. Set $\phi (X):= \cap \psi^{-1}(X)$. Since $\mathcal F$ is closed under intersections, $\overline {\mathcal F}$ is closed under intersections too (Claim \ref{claim:wedge1}). By definition, $\psi$ preserves intersections. It follows that 
$\psi(\phi(X))=X$. From this fact, $\phi(X)$ is the least member $X'$ of $\mathcal H$ such that $\psi (X')=X$. This and the fact that $\psi$ preserves intersections imply that $\phi$ is order preserving.  \end{proof}
\begin{remark} The map  $\phi$ in the proof of Theorem \ref{thm:topologicalbouchet} above does not need to be continuous. For an example, take  $R:=(P, \not \geq , P)$,  $R':= (P', \not \geq, P')$ where $P$ and $P'$ are two posets type $1+\omega^*$ and $(1\oplus 1)+\omega^*$ respectively (here, $1\oplus 1$ denotes a $2$-element antichain) and,  as a coding from $R$ to $R'$,  the  pair $(f,f)$ where $f$ is an embedding from $P$ into $P'$.\end{remark}

From Theorem \ref{thm:topologicalbouchet}, Lemma \ref{lem:ideals2} and Lemma \ref{lem:product}, we derive the following result. 
\begin{proposition} \label {prop:embedproduct} If an incidence structure $R:=(E, \rho, F)$ has a coding in $(Q, \leq, Q)$, where  $Q:=\Pi_{i\in I} C_i$ is a finite product of chains, then $\overline {Gal(R)}$ is embeddable in the product $\Pi_{i\in I} \mathbf {I}(C_i)$.
\end{proposition}
\begin{proof} Set $C'_i:= 1+C_i$ for each $i\in I$ and $Q':=\Pi_{i\in I} C'_i$.  A coding from $(P, \leq , P)$ in $(Q, \leq ,Q)$ induces a coding from $(P, \leq , P)$ in $(Q', \leq ,Q')$. According to Theorem \ref{thm:topologicalbouchet}, such a coding yields an embedding from $\overline {Gal(R)}$ into $\overline {N(Q')}=\overline {Gal((Q',\leq, Q'))}$. The poset $Q'$ is a join-semilattice with a least element, hence according to Lemma \ref{lem:ideals2},  
$\overline {N(Q')}= \mathcal J(Q')$. According to Lemma \ref{lem:product}, $\mathcal J(Q')$ is isomorphic to $\Pi_{i\in I} \mathcal J(C'_i)$. To conclude, observe that $\mathcal J(C'_i)=\mathbf {I}(C_i)$. \end{proof}

We need also the following properties:

\begin{lemma}\label{lem:trivia}Let  $(f,g)$ be a coding from $R:= (E, \rho, F)$ into  $R':= (E', \rho', F')$ and $(\rho'_i)_{\in I}$ such that $\rho':=\bigcap_{i\in I} \rho'_{i}$ then  $\rho= \bigcap_{i\in I} \rho_{i}$ where $\rho_i:=\{(x,y)\in E\times F: f(x)\rho'_i g(y)\}$.
\end{lemma}
The proof is immediate.
\begin{lemma} \label{lem:dimensional} Let $R:= (E, \rho, F)$ be an incidence structure. 
\begin{enumerate}
\item If $\rho= \bigcap_{i\in I} \rho_{i}$ where each $\rho_i$ is an incidence relation from $E$ to $F$, then $Gal(R)$ is embeddable  in $T:= \Pi_{i\in I}Gal(R_i)$ where $R_i:= (E, \rho_i, F)$. 
\item If $Gal(R)$ is embeddable  in a product $C:= \Pi_{i\in I}C_i$ of posets, then $\rho= \bigcap_{i\in I} \rho_{i}$ where each $\rho_i$ is an incidence relation from $E$ to $F$ such that $Gal((E, \rho_i,F))$ is embeddable in $N(C_i)$. 
\end{enumerate}
\end{lemma}
\begin{proof}
$(1)$. Let $A\subseteq E$. Set $\varphi(A):= (L_{R_{i}}\circ U_{R_{i}}(A))_{i\in I}$. Clearly:
\begin{equation}\label{eq:implication}
A\subseteq B\;  \text{implies}\;  \varphi(A)\subseteq \varphi(B).
\end{equation}
Hence $\varphi$ is an order-preserving map from $\mathfrak P(E)$ into $T$. In particular, its restriction to $Gal(R)$ is order-preserving. The fact that this is an embedding is an immediate consequence of the following:  \begin {claim}\label{claim:equality}
\begin{equation} A=  \bigcap_{i\in I} L_{R_{i}}\circ U_{R_{i}}(A) \;\text{provided that}\;  A=L_{R}\circ U_{R}(A).
\end{equation}
\end{claim}
\noindent{\bf Proof of Claim \ref{claim:equality}.}
From $\rho \subseteq \rho_i$ for all $i\in I$, we have $A\subseteq  \bigcap_{i\in I} L_{R_{i}}\circ U_{R_{i}}(A)\subseteq \bigcap_{i\in I} L_{R_{i}}\circ U_{R}(A)$. From $\rho= \bigcap_{i\in I} \rho_{i}$ we get $\bigcap_{i\in I} L_{R_{i}}(B)=L_R(B)$ for every $B\subseteq F$. Applying this to $B:=U_R(A)$, we get $A\subseteq  \bigcap_{i\in I} L_{R_{i}}\circ U_{R_{i}}(A)\subseteq \bigcap_{i\in I} L_{R_{i}}\circ U_{R}(A)=L_R\circ U_R(A)$. The claim follows immediately. \endproof

\noindent $(2)$. Let $c: Gal(R)\rightarrow C$ be an embedding and $p_i: C\rightarrow C_i$ be the $i-th$-projection. Set $f(x):= L_R\circ U_R(x)$ for $x\in E$ and   $g(y):= L_R(y)$ for $y\in F$. Set $f_i := p_i\circ c\circ f$, $g_i := p_i\circ c\circ g$ and  $\rho_{i}:=\{(x,y)\in E\times F: f_i(x)\leq_i g_i(x)\}$. Then $(f_i, g_i)$ is a coding from $R_i:= (E, \rho_i, F)$ into $(C_i, \leq_i, C_i)$. Thus, from Lemma \ref{lem:coding-embedding},  $Gal(R_i)$ is embeddable in $Gal((C_i, \leq_i, C_i))=N(C_i)$. To conclude observe that  $(f,g)$ is a coding from $R$ into $(Gal(R), \subseteq, Gal(R))$,  hence $\rho= \bigcap_{i\in I} \rho_{i}$. \end{proof}

\subsection{Ferrers relations, interval orders and  dimensions}
 Let $R:= (E, \rho,  F)$ be an incidence structure. The binary relation $\rho$ from $E$ to $F$ is a \emph{Ferrers relation} if for every $x,x'\in E$, $y,y'\in F$, $x\rho y$ and $x'\rho y'$ imply $x\rho y'$ or $x'\rho y$. As it is well known, $\rho$ is Ferrers if and only if $Gal(R)$ is a chain. It follows from  Bouchet's theorem that \emph{$Gal(R)$ is a chain if and only if  $R$ has a coding into $(C, \leq, C)$ where $C$ is a chain}.

Let $C$ be a chain, an \emph{interval} of $C$ is any subset $I$ of $C$ such that $x,y\in I, z\in C$ and $x<z<y$ imply $z\in I$. One may order the set $Int (C)$ of non empty intervals of $C$ by setting $I<J$ if  $ x < y$ for all $x\in I$ and $y\in J$. Let $P$ be a poset; the order on $P$ is an \emph{interval order},  and by extension $P$ too,  if $P$ is isomorphic to a subset of  $Int (C)$ for some chain $C$. We recall that:

\begin{lemma}\label{lem:intervalorder}  A poset $P$  is an interval order if and only if $(P,<, P)$ is a Ferrers relation, or equivalently $(P, <, P)$ has a coding into a chain. 
\end{lemma}

   Let $\mathcal F$, resp. $\mathcal J$, be the class of Ferrers relations, resp. interval orders. We recall that the \emph {Ferrers dimension} of  an incidence structure $R:= (E, \rho, F)$  is the least cardinal $\kappa$ such that $\rho$ is the intersection of $\kappa$ \emph{Ferrers relations} from $E$ to $F$. We denote it by $\mathcal F-dim(R)$. The \emph{interval dimension} of $P$ is the smallest cardinal $\kappa$ such that the order on $P$ is the intersection of $\kappa$ interval orders. We denote it by $\mathcal I-dim (P)$. 
  We recall two basic results relating theses notions, due to Bouchet \cite{bouchet} and Cogis \cite{cogis}, namely:
 \begin{equation}\label{ferrerlarge}
\mathcal F-dim((P, \leq, P))=dim (P) 
\end{equation}
and
\begin{equation}\label{ferrerstrict}
\mathcal F-dim ((P, <,P))=\mathcal I-dim(P)
\end{equation}
for every poset $P$.

These three notions of dimension: order dimension, Ferrers dimension and interval dimension are based on three classes of structures: chains, Ferrers relations and interval orders and are expressible in terms of Galois lattices. Replacing theses these classes by others yield other notions of dimension that we discuss at the end of this section. \subsection{ Bipartite posets} 
A poset is \emph{bipartite} if this is the union of two antichains. We recall the following result:
\begin{lemma}
Let $Q$ be a bipartite poset. Then
\begin{equation} \label{bipartite}
\mathcal{I}-dim(Q)\leq dim(Q)\leq \mathcal{I}-dim(Q)+1.
\end{equation}
\end{lemma}

Let $R:= (E, \rho, F)$ be an incidence structure. The \emph{bipartite poset associated to} $R$, denoted by $B(R)$,  is the poset whose base set is $E':= E\times \{0\}\cup F\times\{1\}$ ordered by:
$$(x,i)<(y,j)\;  \text{if } \; (x,y)\in \rho \; \text{ and}  \; i<j.$$

If $P:=(E, \leq)$ we set $B(P):= B(E, \leq, E)$ and $B(\check P):= B(E,<, E)$. The posets $B(P)$ and $B(\check P)$ are respectively called the  \emph{split} and the \emph{open split} of $P$. 

We note that  if $R$ is an incidence structure then $R$ has a coding into $(B(R), \leq , B(R))$ as well as in $(B(R), < , B(R))$. In particular: \begin{equation}\label{eq:edding}
Gal(R) \; \text{is embeddable into}\; N(B(R)).
\end{equation}

As a corollary of (\ref{eq:edding}) it turns out that for every poset $P$:\begin{equation}\label{claim:embedding}
N(P) \; \text{is embeddable into}\;  
N(B(P)).
\end{equation}

Note also that:

\begin{lemma}\label{claim:embedding2}
If $P$ is a poset, $B(P)$ is embeddable in a product $P\times C$
where $C$ is a chain of the form $D+D$, the order type of $D$ being given by any  linear extension of $P^*$.
\end{lemma}
We will use the following easy fact:\begin{lemma}\label{lem:easyfact} Let $R:=(E, \rho, F)$ and $R':= (E', \rho, F')$ be two incidence structures. Every coding $(f,g)$ from $R$ to $R'$ such that $f$ and $g$ are one to one induces an embedding of $B(R)$ in $B(R')$. The converse holds if for every $x\in E$ there is some $y\in F$ such that $(x,y)\in \rho$. 
\end{lemma}

We also recall the following result of Bouchet and Cogis: 
\begin{equation}\label{bouchet-cogis}
\mathcal{F}-dim (R)=\mathcal{I}-dim(B(R))=dim(Gal(R)).
\end{equation}
The first equality in (\ref{bouchet-cogis}) added to equality  (\ref{ferrerlarge}) yields: 
\begin{equation}\label{eq:bip}
dim (P)=\mathcal{I}-dim(B(P)).
\end{equation}
 
 Similarly, the first equality in (\ref{bouchet-cogis}) added to equality  (\ref{ferrerstrict}) yields:
\begin{equation}\label{eq:bipstrict}
\mathcal{I}-dim (P)=\mathcal{I}-dim(B(\check P)).
\end{equation}

Inequalities (\ref{bipartite}) with equality (\ref{eq:bip}) yield 

 \begin{equation} \label{eq:bip1}
dim(P)\leq dim(B(P))\leq dim (P)+1.
 \end{equation}
 
 Similarly, inequalities (\ref{bipartite})  with equality (\ref{eq:bipstrict}) yield
 \begin{equation} \label{eq:split1}
\mathcal{I}-dim(P)\leq dim(B(\check P))\leq \mathcal{I}-dim (P)+1.
 \end{equation}
 
 Inequalities (\ref{eq:bip1}) are due to  Kimble (cf. \cite{trotter-moore}). 
 
Let $\underline {2}.P$ the ordinal product of the two-element chain $\underline 2$ by a poset $P$. This is the set of pairs $(x,i)$, with $x\in P$, $i\in 2$, lexicographically ordered (that is $(x,i)\leq (x',i')$ if either $x<x'$ or $x=x'$ and $i<i'$). 
\begin{lemma}\label{bisplit}
Let  $P$ be a poset and  $Q:= \underline 2.P$ then:
\begin{enumerate}
\item$B(P)$ is embeddable in $B(\check Q)$.
\item $B(\check P)$ is embeddable in $B(Q)$.
\end{enumerate}\end{lemma}
\begin{proof}
Item (1). Let $f$ and $g$ be the maps from $P$ to $Q$ defined by $f(x):=(x,0)$ and $g(x):=(x,1)$. Then $(f,g)$ is a one-to one coding of $(P, \leq , P)$ in $(Q,<,Q)$. This coding induces an embedding from $B(P)$ in $B(\check Q)$. 

Item(2). Let $f'$ and $g'$ be the maps from $P$ to $Q$ defined by $f'(x):=g(x)$ and $g'(x):=f(x)$. Then $(f',g')$ is a one-to one coding of $(P, < , P)$ in $(Q,\leq,Q)$. This coding induces an embedding from $B(\check P)$ in $B(Q)$. 
\end{proof}
\begin{proposition} \label{open-nonopensplit}$B(P)$ and $B(\check P)$ are embeddable in each other whenever  $\underline {2}.P$ is embeddable in $P$.
\end{proposition}
\begin{proof} Let $Q:=\underline 2.P$. Suppose that $Q$ is embeddable in $P$. Then $B(Q)$ is embeddable in $B(P)$. According to item (2) of Lemma \ref{bisplit},  $B(\check P)$ is embeddable in $B(Q)$. Hence $B(\check{P})$ is embeddable in $B(P)$. Similarly, $B(\check Q)$ is embeddable in $B(\check P)$.  According to item (1) of Lemma \ref{bisplit},  $B( P)$   is embeddable in $B(\check Q)$. Hence $B(P)$ is embeddable in $B(\check P)$.\end{proof}

\subsection{A relativisation of the notions  of dimension}
Let $\mathcal R$ be a class of incidence structures and let $R: = (E, \rho, F)$ be an incidence structure.
If $\rho$ is the intersection of incidence relations $\rho_i$ such that  $(E, \rho_i, F)\in \mathcal R$, we define the $\mathcal R$-\emph{dimension} of $R$,  that we denote by $\mathcal R-dim(R)$, as the least cardinal $\kappa$ such that $\rho$ is the intersection of $\kappa$ such relations. Let  $\mathcal D$  be a class of posets and let $P$ be a poset. If the order $\leq$  is the intersection of orders $\leq_i$ such that  $(E, \leq_i)\in \mathcal D$,  the $\mathcal D$-\emph{dimension} of $P$,  that we denote by $\mathcal D-dim(P)$, is the least cardinal $\kappa$ such that $\leq$ is the intersection of $\kappa$ such orders. If the poset  $P$ is embeddable in a product of members of $\mathcal D$ we denote by $\mathcal D-\pi dim(P)$ the least cardinal $\kappa$ such that $P$ is embeddable in a product of $\kappa$ members of $\mathcal D$.  For example,  if $\mathcal R$ is the class $\mathcal F$ of Ferrers relations, $\mathcal R-dim(R)$ is the Ferrers dimension of $R$. If $\mathcal D$ is the class $\mathcal L$ of chains,  $\mathcal D-dim(P)$ is the order dimension of $P$ and if $\mathcal D$ is the  class of interval orders, $\mathcal D-dim(P)$ is the interval dimension of $P$. 

\begin{definition} A class  $\mathcal C$  of posets is \emph{dimensional} if:
 \begin{enumerate} 
 \item $\underline 2 \in \mathcal C$.
 \item If $C\in \mathcal C$ and $C'$ is embeddable in $C$ then $C'\in \mathcal C$.
 \item If $C\in \mathcal C$ then $N(C)\in \mathcal C$. 
 \end{enumerate}
\end{definition}
%\begin{examples} The class of posets, of scattered posets, of chains, of scattered chains are dimensional. 
%\end{examples}

Let $(C_i)_{i\in I}$ be a family of posets such that $I$ is equipped with a well-ordering. The \emph{lexicographical product} of this family is the poset denoted  $\bigodot _{i\in I}C_i$ whose underlying set is the cartesian product $\Pi_{i\in I} C_i$, the ordering being defined by:
$$(x_i)_{i\in I}\leq (y_i)_{i\in I}$$ 
if either $(x_i)_{i\in I}= (y_i)_{i\in I}$ or $x_{i_0}<y_{i_0}$ where $i_0$  is the least $ i \in I$ such that $x_{i_0}\not =y_{i_0}$. 

\begin{proposition} \label{prop:dimension}Let $\mathcal C$ be a dimensional class of posets and 
$Gal^{-1}(\mathcal C)$ be the class of incidence structures $S$ such that $Gal(S)\in \mathcal C$. 
  Then:
\begin{enumerate}[{(i)}]
\item \label{itemi} $Gal^{-1}(\mathcal C)-dim(R) = {\mathcal C}-\pi dim(Gal(R))$ for every  incidence structure $R:= (E, \rho, F)$.  
\item \label{itemii}$Gal^{-1}(\mathcal C)-dim((P, \leq ,P))={\mathcal C}-\pi dim(P)$ for every poset $P$.
\item \label{itemiii} If $\mathcal I ({\mathcal C})$ is the class of posets $(L,\leq)$ such that $Gal((L, <, L))\in \mathcal C$ then ${\mathcal C}-\pi dim(Gal((P, <, P)))\leq {\mathcal I({\mathcal C})}-dim(P)\leq {\mathcal C}-dim(Gal((P, <, P)))$. 
\end{enumerate}
 Let  $\kappa$ be a cardinal. If $\bigodot _{i\in I}C_i\in \mathcal C$ whenever $(C_i)_{i\in I}$ is a family of members of $\mathcal C$ such that $\vert I\vert < \kappa$ then:
 \begin{enumerate}[{(i')}]
 \item \label{itemi'}${\mathcal C}-dim(P)={\mathcal C}-\pi dim(P)$  for every poset $P$ such that   ${\mathcal C}-\pi dim(P)<\kappa$.
 \item \label{itemii'} ${\mathcal I({\mathcal C})}-dim(P)={\mathcal C}-dim(Gal((P, <, P)))$  
for every poset $P$ such that ${\mathcal C}-\pi dim(Gal((P, <,P)))<\kappa$.

\end{enumerate}
\end{proposition}

\begin{proof} Observe that since a poset $P$ is embeddable in the power set $\mathfrak P(P)$ ordered by inclusion, and since this poset is isomorphic to the power $\underline 2^{P}$, $P$ is embeddable in a power of $\underline 2$. Since $\underline 2\in \mathcal C$, ${\mathcal C}-\pi dim (P)$ is well-defined. 

\noindent Item $(i)$. Let $\kappa:={\mathcal C}-\pi dim(Gal(R))$. According to the  observation above, this quantity is well-defined.  Let $C:= \Pi_{i\in I}C_i$ be a product of $\kappa$ members of $\mathcal C$ such that $Gal(R)$ is embeddable in $C$.  According to Item $(2)$ of Lemma \ref{lem:dimensional}, 
$\rho= \bigcap_{i\in I} \rho_{i}$ where each $\rho_i$ is an incidence relation from $E$ to $F$ such that $Gal((E, \rho_i,F))$ is embeddable in $N(C_i)$. Since $\mathcal C$ is dimensional,  $Gal((E, \rho_i,F))\in \mathcal C$, hence  $Gal^{-1}(\mathcal C)-dim(R)$ is well-defined and $Gal^{-1}(\mathcal C)-dim(R) \leq {\mathcal C}-\pi dim(Gal(R))$. The converse inequality follows immediately from Item $(1)$ of Lemma \ref{lem:dimensional}. 

\noindent Item $(ii)$.  We have $N(P)=Gal((P, \leq, P))$. Hence, from Item $(i)$, we have $Gal^{-1}(\mathcal C)-dim((P, \leq ,P))={\mathcal C}-\pi dim(N(P))$. Since $P$ is embeddable in $N(P)$,  ${\mathcal C}-\pi dim(P)\leq {\mathcal C}-\pi dim(N(P))$. To get the converse inequality, note that if $P$ is embeddable in a product $C:=\Pi_{i\in I}C_i$ then, since each $C_i$ is embeddable in $N(C_i)$, $C$ is embeddable in  $C':= \Pi_{i\in I}N(C_i)$, hence $P$ is embeddable in $C'$. Since $C'$ is a complete lattice, $N(P)$ is embeddable in $C'$.  From the  fact that $\mathcal C$ is dimensional, $C'\in \mathcal C$. The result follows. 

\noindent Item $(iii)$.  Set $R:= (P, <,P)$. We  prove first  the second inequality. 
\begin{claim} \label {claim:interval}$\mathcal I({\mathcal C})-dim(P)\leq {\mathcal C}-dim(Gal(R))$. 
\end{claim}
\noindent{\bf Proof of Claim \ref {claim:interval}.} Let $(f,g)$ be  the coding from $R$ into $(Gal(R), \subseteq, Gal(R))$ defined by  $f(x):= L_R\circ U_R(x)$ for $x\in P$ and   $g(y):= L_R(y)$ for $y\in P$. We have:
\begin{equation} \label{eq:intorder} g(x)\subset f(x)
\end{equation}
 for all $x\in P$. Indeed, since $x\not <x$, we have $f(x)\not \subseteq g(x)$; on an other hand we have $g(x)=L_R(x)\subseteq L_R\circ U_R(x)=f(x)$. Now, let  $\mathcal L'$  be an order extending the inclusion order on $Gal(R)$. Set  $\mathcal L:= \{(x,y)\in P: (f(x),g(y))\in \mathcal L'\}$. Then $\mathcal L$  is   irreflexive and transitive. Indeed, according to (\ref{eq:intorder}) we have $g(x)\subset f(x)$,thus $(g(x), f(x))\in \mathcal L'$. This implies that $(f(x),g(x))\not \in \mathcal L$, hence $(x,x)\not \in \mathcal L$, proving that $\mathcal L$ is irreflexive. Let $(x,y), (y,z)\in \mathcal L$. Since $g(y)\subset  f(y)$, $(g(y), f(y))\in \mathcal L'$. This easily yields that $(x, y)\in \mathcal L$, thus $\mathcal L$ is transitive. If, moreover $(Gal(R), \mathcal L')\in \mathcal C$, $Gal((P, \mathcal L, P))\in \mathcal C$. Indeed, $(f,g)$ is a coding from  $(P, \mathcal L, P)$ into $(Gal(R), \mathcal L', Gal(R))$. Hence, from Bouchet's theorem (cf. Corollary \ref{lem:coding-embedding}),  $Gal((P, \mathcal L, P))$ is embeddable into $Gal((Gal(R), \mathcal L', Gal(R)))= N((Gal(R), \mathcal L'))$. Since $\mathcal C$ is dimensional, if $(Gal(R), \mathcal L')\in \mathcal C$, $N((Gal(R), \mathcal L'))\in \mathcal C$ too. Thus $Gal((P, \mathcal L, P))\in\mathcal C$. With that,   our claim follows from Lemma \ref{lem:trivia}. \endproof

\begin{claim}\label{claim:idim}
$Gal^{-1}(\mathcal C)-dim(R)\leq \mathcal I({\mathcal C})-dim(P)$.
\end{claim}
\noindent{\bf Proof of Claim \ref{claim:idim}.}
Trivial.
\endproof. 

From $(i)$ we have $Gal^{-1}(\mathcal C)-dim(R) = {\mathcal C}-\pi dim(Gal(R))$. Thus, with Claim \ref{claim:idim}, ${\mathcal C}-\pi dim(Gal(R))\leq \mathcal I({\mathcal C})-dim(P)$. This is the first inequality. With that, the proof of Item $(iii)$ is complete. 

\noindent Item $(i')$. We have ${\mathcal C}-\pi dim(P)\leq {\mathcal C}-dim(P)$ without any condition on $\mathcal C$. Indeed, if the order $\leq$ on $P$ is the intersection of a family $(\leq_i)_{i\in I}$ of orders on $P$, the map $\delta: P\rightarrow P^I$ defined by  $\delta(x)(i):=x$  is an embedding of $P$  in the direct product 
$\Pi_{i\in I} P_i$ where $P_i:=(P, \leq_i)$. Conversely, suppose that there is an embedding from $P$ in a direct product $Q:= \Pi_{i\in I} P_i$, with $P_i\in \mathcal C$. Let $P'$ be the image of $P$.
\begin{claim}\label{claim:lexico}The order on $Q$ is the intersection of $
\vert I\vert $ orders $\leq_i$ such that $(Q, \leq_i)\in \mathcal C$. 
\end{claim}
\noindent{\bf Proof of Claim \ref{claim:lexico}.}
For each $i\in I$, choose a well-ordering $\mathcal L_i$ on $I$  for which $i$ is the first element and let $Q_i$ be the lexicographical product of the $P_i$'s indexed by $L_i:= (I, \mathcal L_i)$. The order on $Q$ is the intersection of the orders of the $Q_i$'s. If each $P_i$ belong to $\mathcal C$, then with our hypothese  on $\mathcal C$, the $Q_i$'s belong to $\mathcal  C$. \endproof

Now, the order on $P'$  is the intersection of the orders induced on $P'$ by the $Q_i$'s. Since $\mathcal C$ is dimensional, these orders belong to $\mathcal C$, hence  ${\mathcal C}-dim(P)\leq \vert I\vert$.
Thus Item $(i')$ holds.

Item $(ii')$.  Apply Item $(i')$ to $Gal((P, <,  P))$ and use Item $(iii)$.

With this, the proof of Proposition \ref{prop:dimension} is complete.
\end{proof}

 Since the class of chains is preserved under lexicographical products, Proposition \ref{prop:dimension} applied to $\mathcal C:= \mathcal L$ yields formulas  (\ref{ferrerlarge}) and (\ref{ferrerstrict}).

\section{Scattered posets and scattered topological spaces}\label{section:scattered}

A poset $P$, or its order as well,  is \emph{scattered} if it does not contain a subset ordered as the chain $\eta$ of rational numbers; in other words, the chain $\eta$   is not embeddable in $P$.  A topological space is {\it scattered} 
if every non-empty subset has at least an isolated point(w.r.t. the induced topology). Sometimes, to avoid confusion, we use the terms \emph{ordered scattered} and  \emph{topologically scattered}. These  two notions are quite
related. This is particularly the case when the order and the 
topology are defined on the same universe. For an example, if the ordering 
is linear and the topology is the interval-topology, the chain is 
complete if and only if  the space compact (Hausdorff). Moreover, \emph{if $C$ is  a complete  chain, the conditions that  $C$  is order-scattered, 
$C$ is topologically scattered, $C$ is order isomorphic to ${\bf I } (D)$, where $D$ is a scattered chain, are equivalent}. From this fact follows that \emph{a chain $D$ is order scattered if and only if its MacNeille completion $N(D)$ is order scattered}. 

The class $\mathcal S$ of scattered posets is closed downward, that is if $P\in \mathcal S$ and $Q$ is embeddable in $P$ then $Q\in \mathcal S$. Furthermore, it is closed under finite direct product and under finite lexicographical product. In particular, the class $\mathcal L_{\mathcal S}$ of scattered chains is preserved under finite lexicographical product. This property, and $(i')$ of Proposition \ref{prop:dimension}, yield an important known fact: 

\begin{proposition}\label{cor:product}  Let $n$ be a positive integer. An ordered set $P$ is embeddable in a product of $n$ scattered
chains if and only if the order on $P$ is the intersection of $n$ scattered linear orders.
\end{proposition}
%Let $n$ be a non negative  integer; let  $\mathcal {L}(n)$, resp. $\mathcal {L}_{\mathcal S}(n)$,  be the class of posets $P$ whose order is the the  intersection of at most $n$ linear orders, resp. at  most $n$ scattered linear orders (hence $\mathcal {L}(0)=\mathcal {L}_{\mathcal S}(0)$ is made of the empty chain and the one element chains, whereas $\mathcal {L}(1)=\mathcal {L}$ , resp. $\mathcal {L}_{\mathcal S}(1)=\mathcal {L}_{\mathcal S}$).

\subsection{Scattered dimensions}
Let $\mathcal F_{\mathcal S}$, resp. $\mathcal I_{\mathcal S}$, be the class of incidence structure $R$, resp. posets $P$, such that the Galois lattice $Gal(R)$, resp.  $Gal((P,<,P))$ belongs to $\mathcal L_{\mathcal S}$. 

The following lemma completes the  analogy between $\mathcal I_{\mathcal S}$ and $\mathcal I$
\begin{lemma}\label{lem:intervalorder}  A poset $P$ belongs to $\mathcal I_{\mathcal S}$ if and only if $P$ is isomorphic to a subset of  $Int (C)$ for some scattered chain $C$. 
\end{lemma}
%\begin{proof}Let Say that two elements $x$ and $y$ of $P$ are equivalent if they have the same predecessors and the same successors. 
%\end{proof}
%Let $\kappa$ be a cardinal.   If $R:=(E, \rho, F)$ is an incidence structure, resp. if $P$ is a poset, we set  $\mathcal F_{\mathcal S}-dim(R):=\kappa$, resp.  $\mathcal I_{\mathcal S}-dim(P):=\kappa$,  if $\rho$ is the intersection of $\kappa$ relations $\rho'$ such that $(E, \rho',F)\in \mathcal F_{\mathcal S}$, resp. the order is the intersection of $\kappa$  orders  $\leq'$ such that $(P, \leq')$ belongs to $\mathcal I_{\mathcal S}$  and no smaller cardinal has this property.  

Let $n$ be an integer, we denote by $\mathcal F(n)$, resp. $\mathcal I(n)$, resp. $\mathcal L(n)$ the class of incidence structures $R$, resp. of posets $P$ such that $\mathcal F-dim (R)\leq n$, resp. $\mathcal {I}-dim (P)\leq n$, resp. $dim(P)\leq n$.  We define $\mathcal F_{\mathcal S}(n)$, resp. $\mathcal I_{\mathcal S} (n)$, resp. $\mathcal L_{\mathcal S}(n)$, accordingly. 

\begin{theorem} \label{thm:extension}Let $n$ be an  integer and let $R$ be an incidence structure, resp. a poset $P$. Then $R\in \mathcal F_{\mathcal S}(n)$, resp. $P\in \mathcal I_{\mathcal S}(n)$, resp.   $P\in \mathcal L_{\mathcal S}(n)$, if and only if $Gal(R)$, resp. $Gal((P, <, P))$, resp. $N(P)$, belongs to $\mathcal L_{\mathcal S}(n)$.   
\end{theorem}
\begin{proof} We apply Proposition  \ref{prop:dimension} with $\mathcal C:=\mathcal L_{\mathcal S}$.  Since $\mathcal F_{\mathcal S}=Gal^{-1}(\mathcal L_{\mathcal S})$,   Item (\ref{itemi}) yields $\mathcal F_{\mathcal S}-dim(R) = \mathcal L_{\mathcal S}-\pi dim(Gal(R))$ for every  incidence structure $R:= (E, \rho, F)$. From Proposition \ref{cor:product}, ${\mathcal L_{\mathcal S}}-dim(Gal(R))={\mathcal L_{\mathcal S}}-\pi dim(Gal(R))$ provided that   ${\mathcal L_{\mathcal S}}-\pi dim(Gal(R))<\omega$. Thus if $R\in \mathcal F_{\mathcal S}(n)$,  ${\mathcal L_{\mathcal S}}-dim(Gal(R))\leq n$, that is $Gal(R)\in\mathcal L_{\mathcal S}(n)$. The converse follows from the fact that ${\mathcal L_{\mathcal S}}-\pi dim(Gal(R))\leq {\mathcal L_{\mathcal S}}-dim(Gal(R))$. Set $R:= (P, <,P)$. Since $\mathcal I_{\mathcal S}=\mathcal I(\mathcal C)$, Item (\ref{itemiii}) yields ${\mathcal L_{\mathcal S}}-\pi dim(Gal(R)\leq \mathcal I_{\mathcal S}-dim(P)\leq {\mathcal L_{\mathcal S}}-dim(Gal(R)$. Thus, if    $P\in \mathcal I_{\mathcal S}(n)$,  ${\mathcal L_{\mathcal S}}-\pi dim(Gal(R)\leq n$. Since ${\mathcal L_{\mathcal S}}-\pi dim(Gal(R)={\mathcal L_{\mathcal S}}-dim(Gal(R)$,  $Gal(R)\in\mathcal L_{\mathcal S}(n)$. Again, the converse follows from the fact that ${\mathcal L_{\mathcal S}}-\pi dim(Gal(R))\leq {\mathcal L_{\mathcal S}}-dim(Gal(R))$. Now, set $R:=(P,\leq, P)$. Combining Item $(i)$ and Item $(ii)$,  we get ${\mathcal C}-\pi dim(Gal(R))=Gal^{-1}(\mathcal C)-dim(R)={\mathcal C}-\pi dim(P)$, hence ${\mathcal C}-\pi dim(N(P))={\mathcal C}-\pi dim(P)$. Since ${\mathcal C}-\pi dim(P)\leq {\mathcal C}-dim(P)$,  if $P\in \mathcal L_{\mathcal S}(n)$, ${\mathcal C}-\pi dim(N(P))\leq n$. With Item $(i')$ we get $N(P)\in \mathcal L_{\mathcal S}(n)$. The converse is similar. 
\end{proof}

 With these notations, one may ask: 
 \begin{questions}  Let $P$ be a poset and $R$ be an incidence structure. 
 \noindent \begin{enumerate}[{(i)}]

\item If $\mathcal L_{\mathcal S}-dim(P)$ is finite does  ${\mathcal L_{\mathcal S}}-dim(P)= dim(P)$?

\item If  $\mathcal I_{\mathcal S}-dim(P)$ is finite does $\mathcal I_{\mathcal S}-dim(P)= \mathcal{I}-dim(P)$? 

\item If  $\mathcal F_{\mathcal S}-dim(R)$ is finite does $\mathcal F_{\mathcal S}-dim(R)= \mathcal{F}-dim(R)$? 
\end{enumerate}
\end{questions}

Question $(i)$ is just a reformulation of Question \ref{que:scatt}. 

With the help of Theorem \ref{thm:extension}, one can show that a positive answer to $(i)$ is equivalent to a positive answer to $(iii)$ and implies a positive answer to $(ii)$. 

\subsection{Topologically scattered spaces and Galois lattices}
. 
 \begin{lemma}\label{toposcatt}
 \begin{enumerate}
 \item \label{rudin}The continuous image of a compact scattered space is scattered. 
 \item A finite product of scattered topological spaces is scattered.
 \item A Priestley space which is topologically scattered is order scattered. 
 \end{enumerate}
 \end{lemma}
 
 The first fact is non-trivial, it is due to W.Rudin.  The second and third fact are easy and well-known.

\noindent{\bf Remark} If $L$ is a topologically scattered algebraic lattice, the algebraic lattice $\overline {N(L)}$ is not necessarily topologically scattered.  A  topologically scattered algebraic lattice $L$ containing an infinite independent set $X$ will do. Indeed, recall that a subset $X$ of a lattice $L$ is \emph{independent} if $x\not \leq \vee X$ for every $x\in X$, $F\in [X\setminus \{x\} ]^{<\omega}$. Furthermore, if $X$ is independent,  $\mathfrak P(X)$ is embeddable in $\mathcal J(L)$. Thus, if $X$ is infinite,  $\mathcal J(L)$ is not order scattered. Since  $\overline {N(L)}=\mathcal J(L)$, this set is not topologically scattered. For that, let $P$ be a  countable well-founded poset with no infinite antichain. Set $L= {\bf I} (P)$. Then $L$ is countable, thus topologically scattered.  Since $L$ is distributive,  antichains of join-irreducible members of $L$ are independent subsets of $L$. To get $L$ containing an infinite antichain of join-irreducibles,  take for $P$ the poset made of  $\{ (m,n) \in \N^2: m<n\}$ ordered by setting \begin{equation}\label{eqrado}(m,n)\leq_P (m',n')\; \mbox {if either} \; m=m'\;  \mbox {and}\;  n \leq n'  \; \mbox {or}\;  n<m'
\end{equation} This poset was discovered by R. Rado \cite{rado}.

\begin{lemma}\label{lem:cor1bouch}
Let  $R:= (E, \rho, F)$ and $R':= (E', \rho', F')$ be two incidence structures. If $R$ has a coding into $R'$ and  $\overline{ Gal(R')}$ is topologically scattered then $\overline{ Gal(R)}$ is topologically scattered. \end{lemma}
\begin{proof}
According to Theorem \ref{thm:topologicalbouchet}, $\overline {Gal(R)}$ is the continuous image of a closed subspace of $\overline {Gal(R')}$. From Rudin's result ((\ref{rudin}) of Lemma \ref{toposcatt}) it is topologically scattered. 
\end{proof}

\begin{lemma}\label{lem:pdscatteredchains}
Let  $n$ be an integer and $R:= (E, \rho, F)$ be an incidence structure. If  $Gal(R)$ is embeddable into a  product of  $n$ scattered chains, then $\overline {Gal(R)}$ too. Moreover, $\overline {Gal(R)}$ is topologically scattered.  \end{lemma}
\begin{proof}
Suppose that  $Gal(R)$ is embeddable in  $Q:=\Pi_{i\in I} C_{i}$ with $\vert I \vert =n$. According to Bouchet's theorem (Theorem \ref{bouchet}), $R$ has a coding into $Q$.   According to Theorem \ref{thm:topologicalbouchet}, $\overline{Gal(R)}$ is embeddable in $\overline {Gal((Q, \leq , Q))}=\overline {N(Q)}$. Since $Gal(R)$ has a least element, we may suppose w.l.o.g that $Q$ has a least element, that is each $C_i$ has a least element. Then  $Q$ is a join-semilattice with a least element, hence from Lemma \ref{lem:ideals2}, $\overline {N(Q)}=\mathcal J(Q)$. Since $I$ is finite,  Lemma \ref{lem:product} ensures that   $\mathcal J(Q)$ is order isomorphic to $\Pi_{i\in I} \mathcal J(C_{i})$, that is to  $\Pi_{i\in I} \mathbf {I}(C'_{i})$, where each $C'_i$ is such that $C_i=1+C'_{i}$.  Since the $C'_i$'s are order scattered, the $\mathbf {I}(C'_i)$'s are  order scattered. Their product, being finite,  is scattered too.   This proves the first assertion. The $\mathbf {I}(C'_i)$'s are  in fact topologically scattered. Hence, as  a finite product of scattered spaces, $\mathcal J(Q)$ is  topologically scattered. The second part of the assertion follows from  Lemma \ref{lem:cor1bouch}.
\end{proof}

\section{Proofs of Theorems 1, 2 and 3}\label{section:proofs}

\subsection{Proof of Theorem \ref{thm:second}}
$(i)\Rightarrow (ii)$. Suppose that $(i)$ holds. Then, according to Proposition \ref{cor:product},  $P$ is embeddable in a product $Q:=\Pi_{i\in I}C_i$ of $n$ scattered chains.  In particular  $(P\leq P)$ has a coding into $(Q, \leq Q)$. According to Proposition  \ref{prop:embedproduct},   $\overline {N(P)}$ is embeddable in $\Pi_{i\in I} {\bf I}(C_{i})$. Thus $(ii)$ holds. Moreover, from Lemma \ref{lem:pdscatteredchains}, $\overline {N(P)}$ is topologically scattered. 
$(ii)\Rightarrow (i)$. Suppose that $(ii)$ holds. Since $P$ is embeddable in $\overline {N(P)}$, it is embeddable in a product of $n$ scattered chains. According to Proposition \ref{cor:product}, $(i)$ holds.

\subsection{Proof of Theorem \ref{bip2}}
Suppose that $B(P)\in \mathcal L_{\mathcal S}(n)$. From Theorem  \ref{thm:second}, $N(B(P))\in \mathcal L_{\mathcal S}(n)$. Since from (\ref{claim:embedding}),  $N(P)$ is embeddable into $N(B(P))$, $N(P)\in  \mathcal L_{\mathcal S}(n)$. Hence, $P\in  \mathcal L_{\mathcal S}(n)$. 
Conversely, suppose that $P\in \mathcal L_{\mathcal S}(n)$. In this case, we apply Lemma \ref{claim:embedding2} with   $D\in \mathcal L_{\mathcal S}$. It turns out that  $B(P)\in  \mathcal L_{\mathcal S}(n+1)$. \endproof

\subsection{Proof of Theorem \ref{thm:main}}We prove that if $P$ is one of the ten posets listed in Theorem \ref{thm:main}, either $\overline {N(P)}$ or $\overline {N(P^*)}$ is not topologically scattered. According to Theorem \ref{thm:second} the order on $P$ cannot be the intersection of finitely many scattered linear orders and thus any poset containing a copy of $P$ has the same property. 

Since  for each member $P$ of our list, $P^*$ belongs to our list,  it suffices to check that $\overline {N(P)}$ is not topologically scattered  in the following cases. 

\noindent{\bf Case 1.}  $P\in \{\eta, T_2,\Omega (\eta)\}$.  If $P=\eta$, $\overline {N(P)}={\bf I}(\eta)$.  Topologically, this space is the Cantor set; it is not topologically scattered.  If $P=T_2$, then $\overline {N(P)}$ is made of the binary tree plus the maximal branches of the binary tree and a top element added. These maximal branches form a Cantor space, hence   $\overline {N(P)}$ is not topologically scattered (a strengthening of this fact will be given in Proposition \ref{lem:separdyadic}). If $P= \Omega (\eta)$, the pictorial representation of $\Omega (\eta )$ given in Figure \ref{Omega} show that  $\Omega (\eta )$ is a $2$-dimensional poset, in fact the intersection of a linear  order of type $\omega$ and of a linear order of type $\omega.\eta$. Moreover, as it is easy to see,
${\bf I}(\eta)$ is embeddable in $\mathcal J(\Omega (\eta ))$.
Since $\mathcal J(\Omega (\eta ))\subseteq \overline {down (\Omega (\eta ))}\subseteq \overline {N(\mathcal J(\Omega (\eta ))}$, it follows that $\overline {N(\mathcal J(\Omega (\eta ))}$ is not order scattered, hence not topologically scattered. 

{\bf Case 2. } $P:= B(\check Q)$ where $Q\in  \{\eta, T_2,\Omega (\eta)\}$. 
We deal with the three cases at once. Since   $ {N(Q)}\setminus \overline {Gal((Q, <,Q))}$ is made of isolated points, it follows from Case 1 that  $\overline {Gal((Q, <,Q))}$ is not topologically scattered. Since $(Q, <, Q)$ has a coding into $B((Q, <, Q))=B(\check Q)$, Theorem \ref{thm:topologicalbouchet} yields that  $\overline {Gal( (Q, <, Q))}$ is the continuous image of  $\overline {N(B(\check{Q}))}$. From Rudin's result ((\ref{rudin}) of Lemma \ref{toposcatt}) this latter set cannot be  topologically scattered. \endproof

\begin{lemma} \label{lem:comp}If $P\in\{\eta, T_2,\Omega (\eta)\}$, $B(P)$ and $B(\check P)$ are embeddable in each other. \end{lemma} 
\begin{proof}
Observe that $\underline {2}.P$ is embeddable in $P$ and apply Proposition \ref{open-nonopensplit}.
\end{proof}
\begin{lemma}\label{lem:dimten}
The ten posets listed in Theorem \ref{thm:main} have dimension at most $3$.
\end{lemma}
\begin{proof} Trivially $\eta$ has dimension $1$. As a tree,  $T_2$ has dimension $2$.  Figure \ref{Omega} shows that  $\Omega(\eta)$ has dimension $2$. The poset  $B(\check \eta)$ is defined as the strict product of the chain of rational numbers and the $2$-element chain on $\{0,1\}$ with $0<1$. Hence, this is a $2$-dimensional poset. Let $P\in\{T_2,\Omega (\eta)\}$. Since $P$ has dimension $2$, it follows from equation (\ref{eq:bip1}) that $B(P)$ has  dimension at most $3$. According to  Lemma  \ref{lem:comp}, $B(\check P)$ is embeddable in $B(P)$, thus $B(\check P)$ has  dimension at most $3$. Let  $A:= \{0\}\cup 3\times 2$ ordered so that $0$ is the least element and $(i,j)<(i',j')$ if $i=i' $ and $j<j'$. This poset is a tree obtained by taking the direct sum of three copies of a $2$-element chain and adding a least element. This tree is obviously embeddable in  $T_2$. Every $2$-dimensional poset is embeddable in $\Omega (\eta)$ thus $A$ is also embeddable in $\Omega (\eta)$. Let $X:=\{(0,0)\}\cup \{((i,j),j): i<3, j<2\}$ and $B:=B(\check A)_{\restriction X}$. The poset $B$ is  $3$-dimensional poset (in fact a $3$-irreducible poset).  Since $A$ is embeddable in $P$, $B(\check A)$ is embeddable in $B(\check P)$. Thus, $B(\check P)$ has dimension $3$.  With the fact that  a poset and its dual have the same dimension, our proof is complete. \end{proof}
\begin{lemma} \label{lem:incomp}The ten posets listed in Theorem \ref{thm:main} are pairwise incomparable with respect to embeddability. 
\end{lemma} 
\begin{proof}
Let $X_0:=\eta$, $X_1:=T_2$,  $X_2:=\Omega (\eta)$, $X_3:=B(\check \eta)$, $X_4:=B(\check T_2)$,  $X_5:=B(\check \Omega (\eta))$, $X_6:=(T_2)^*$,  $X_7:=(\Omega (\eta))^*$, $X_8:=(B(\check T_2))^*$,  $X_9:=(B(\check \Omega (\eta)))^*$. We need to prove that $X_i$ is not embeddable in $X_j$ for all pairs $(i,j)$ of distinct elements. Clearly, it suffices to consider the pairs for which $i\leq 5$ and $j\leq 9$. We consider only pairs $(i, j)$ for which a significant argument is needed. For the pair $(1, 2)$ note that $\mathcal J^{\neg \downarrow\hskip -2pt }(X_1)$ is an antichain whereas $\mathcal J^{\neg \downarrow\hskip -2pt }(X_2)$ is a chain. For pairs $(3, j)$, with $j\not \in \{0, 6,8, 9 \}$,  note that $\mathbf I( X_3)$ contains principal initial segments which are infinite whereas $\mathbf I(X_j)$ contains  none.
For  pairs $(i,j)$ such that 
$i\in \{4,5\}$  and $j\not \in \{4,5, 8,9\}$ note that $dim (X_i)=3$ and $dim (X_j)\leq 2$ (Lemma \ref{lem:dimten}). For pairs $(i, j)$ such that 
$i\in \{4, 5\}$, $j\in \{3, 4, 5, 8\}$, then we may write $X_i=B(\check{Y_i})$ and $X_j=B(\check{Y_j})$. If $X_i$ is embeddable in $X_j$, it follows from Lemma \ref{lem:comp} that   $B({Y_i})$ is embeddable in $B(Y_j)$. From Lemma \ref {lem:easyfact} there is a coding from $(Y_i, \leq, Y_i)$ in $(Y_j, \leq, Y_j)$, from which follows that  $\overline{N(Y_i)}$ is embeddable in $\overline{N(Y_j)}$. Since $Y_i=X_{i'}$ for some $i'\leq 2$, this yields that   $\overline{N(X_{i'})}$ is embeddable in $\overline{N(Y_j)}$. Except for the pair $(5,3)$ (which has been previously ruled out) this is clearly impossible. With this last argument, the proof is complete.\end{proof}

\section{Scattered distributive lattices}\label {section:proof4}
In this section, we consider \emph{bounded} distributive lattices, that is distributive lattices with a least and a largest element denoted respectively by $0$ and $1$. If $T$ is a such a lattice, an ideal $I$ is \emph{prime} if its complement $T\setminus I$ is a filter. The \emph{spectrum} of $T$, that we denote $Spec (T)$, is  the subset of $\mathfrak P(T)$ made  of prime ideals of $T$. W.r.t. the topology on $\mathfrak P(T)$, this  a closed subspace of $\mathfrak P(T)$, and with the inclusion order added, this is a Priestley space. The set of order preserving and continuous maps  from $Spec(T)$ onto the two element chain $\underline 2$ is a distributive lattice isomorphic to $T$. This fact is the essence of Priestley duality. We give below the facts we need in order to prove Theorem \ref {thm:distributivelattice}. We only give the proofs or an hint when needed.
The first one is obvious:
\begin{lemma}\label{lem:examplepriestley} If $D$ is a chain with  least and largest elements $0$ and $1$, then as a Priestley space, $Spec(D)$ is isomorphic to ${\bf I}(C)$ where $C:=D\setminus\{0, 1\}$. 
\end{lemma}

\begin{lemma}\label{lem:maximalpriestley} Let $T$ be a distributive lattice and $C$ be a maximal chain  of $Spec(T)$. Then, as a Priestley space,  $C$ is isomorphic to $Spec(D)$ where $D$  a chain, quotient of $T$. 
\end{lemma}
For a proof, note that  the  spectrum of $T$, $Spect(T)$,  is closed under unions and intersections of non-empty chains. Hence $C$ is a complete chain.

We recall that the \emph{width} of a poset $P$, denoted by $width(P)$, is the supremum of the  cardinalities of the antichain of $P$. The following result is due to Dilworth \cite{dilworth}.  
\begin{theorem}\label {lem:spectdim}Let $T$ be a distributive lattice and $n$ be an integer. Then $dim(T)\leq n$ if and only if $width (Spec (T))\leq n$.\end{theorem}

\begin{lemma}Ê\label{lem:spec} Let $T$ be a distributive lattice,  two elements $x, y$ of $T$ such that $x<y$ and $T':= [x, y]$. Then $Spec(T')$ is isomorphic as a Priestley space to $A:=\{J\in Spect(T): x\in J\; \text{and}\; y\not \in J\}$. \end{lemma}
\begin{proof}Let $\phi:A \rightarrow Spec(T')$ and $\theta: Spec(T')\rightarrow A$ defined by setting $\phi(I):=I\cap T'$ and $\theta(I'):=\downarrow I'$ are order preserving, continuous and inverse of each other.
\end{proof}

\begin{lemma}\label{lem:keyscatered}Let $T$ be a distributive lattice. If $Spec(T)$ is not topologically scattered, there is some element $x\in T\setminus\{0,1\}$ such that the spectra of $T'':=\downarrow x$ and $T':=\uparrow  x$ are not topologically scattered. 
\end{lemma}
\begin{proof} Since $Spec(T)$ is not topologically scattered, it contains a perfect subspace. Let $P$ be such a subspace. Since $\vert P\vert \geq 2$, we may pick $J',J''\in P$ such that $J'\not\subseteq J''$. Let $x\in J'\setminus J''$. Then $Spec(T'')$ and $Spec( T')$ are not scattered. Indeed, note first that according to Lemma \ref{lem:spec},  $Spec(T'')= \{J\in Spec (T): x\not \in J\}$ and $Spec( T')= \{J\in Spec (T): x \in J\}$.  Next, observe that the sets   $F':=\{J\in P: x\in J\}$ and $F'':=\{J\in P: x\not \in J\}$ are perfect. Since there are respectively contained in $Spec(T')$ and $Spec(T'')$ the conclusion follows. \end{proof}

\begin{theorem}\label {lem:spectscattered} Let $T$ be a distributive lattice. Then $T$ is order-scattered if and only if $Spec(T)$ is topologically scattered.
\end{theorem}
\begin{proof}  Suppose that $Spec(T)$ is not topologically scattered.  For each pair of elements $x,y$ in $T$ such that $Spec([x,y])$ is not topologically scattered,  Lemma \ref{lem:keyscatered} yields some $z\in ]x, y[$ such that neither $Spec([x,z])$ nor $Spec([z,y])$ is topologically scattered. This fact allows to define an embedding $\phi$ from the set $D$ of dyadic numbers of the $[0,1]$ interval of the real line.   Since $D:= \{\frac{m}{2^n}: n \leq m\in \N\}$ is dense, $T$ is not order scattered.  Conversely, if $T$ is not order scattered,select a non scattered chain and extend it to a maximal chain, say $D$. The natural embedding from $D$ into $T$ yields a continuous surjective map from $Spec(T)$ onto $Spec (D)$. As a Priestley space, $Spec (D)$ is isomorphic to ${\bf I}(C)$ where $C:=D\setminus \{0, 1\}$ (Lemma \ref{lem:examplepriestley}). Since $C$ is not order scattered, $Spec(D)$   is not topologically scattered. According to Rudin's result ((\ref{rudin}) of Lemma \ref{toposcatt}), $Spec(T)$ is not topologically scattered. 
\end{proof}

\subsection{Proof of Theorem \ref{thm:distributivelattice}.} We prove the result for bounded lattices. If $T$ is not bounded, we add a least and a largest element, and apply the result to the resulting lattice. 
$(ii)\Rightarrow (i)$ Apply Proposition \ref{cor:product}.\\
$(i)\Rightarrow (iii)$ Trivial. \\
For the proof of $(iii) \Rightarrow (ii)$, we introduce the following property:\\
\emph{\noindent $(iv)$ $Spec(T)$ is order scattered and $width(Spec(T))\leq n$}. 

We prove successively $(iii)\Rightarrow (iv)$ and $(iv)\Rightarrow (ii)$.\\
$(iii)\Rightarrow (iv)$. Suppose that $(iii)$ holds. Since $T$ is order scattered, Theorem \ref{lem:spectscattered} ensures that  $Spec(T)$ is topologically scattered. With the inclusion order and the topology,  $Spec(T)$ is a Priestley space, hence it is order scattered. Since $dim(T)\leq n$,  Theorem  \ref {lem:spectdim} ensures that  $width (Spec (T))\leq n$. Thus, $(iv)$ holds. \\
$(iv)\Rightarrow (ii)$. Suppose that $(iv)$ holds. Cover $Spec(T)$ with $m$ chains, where $m:=width (Spec (T))$. Extend each of these  chains to a maximal chain of $Spec(T)$.  According to Lemma\ref{lem:maximalpriestley}, each maximal chain $C_i$ is of the form $Spec (D_i)$ where $D_i$ is a chain. Since $Spec(T)$ is order scattered, $C_i$ and hence $D_i$ is order scattered. Let $C:= \oplus_{i<m} Spec (D_i)$ and  $f: C\rightarrow Spec(T)$ defined by setting $f(x,i):=x$. The duality between distributive lattices and their Priestley spaces, yields a lattice  embedding from $T$ into $\Pi_{i<n} D_i$. Hence $(ii)$ holds. \endproof

\section{ Two-dimensional scattered posets}\label{section:two-dimensional}
A linear extension $L$ of an ordered set $P$ is called \emph{separating} if there
are elements $x,y,z\in P$
with $x<_P z$, $y$ incomparable with both $x$ and $z$ but
$x<_L y<_L z$. Let $P$ be an ordered set. If the order of $P$  is the intersection  of two non-separating  linear extensions $C$ and $C'$ of $P$, $C'$ is called  a \emph{complement} of $C$. 
 Dushnik and Miller\cite{dushnik-miller} gave the following
characterization of ordered sets of dimension at most $2$.
\begin{theorem} \label{thm:dim2}Let $P$ be an ordered set, the following properties are
equivalent:
\begin{enumerate}[{(i)}]
\item $dim(P) \leq 2$.
\item There is a linear extension of $P$ which is non-separating.
\item $P$ is embeddable in the family of intervals of some
chain, these intervals being  ordered by inclusion.
\item
The incomparability graph of $P$ is a comparability graph.
\end{enumerate}
\end{theorem}
We mention the following property.
\begin{lemma} \label{basictool} Let $P$ be a  poset of dimension $2$,
$L$ be a non-separating extension of $P$ and $I$ be an initial
segment of $L$. For $z \in P$, we define
$D(z):=(\downarrow z) \cap I$. Then for every
$x,y \in P \setminus I$, with $x \parallel _P y$ we have:
\begin{enumerate}
\item[ 1)] $D(x) \subseteq D(y)$ or $D(y) \subseteq D(x)$,
\item[ 2)] If $D(x) \subset D(y)$ then $y<_L x$.
\end{enumerate}
\end{lemma}
\begin{proof}
1) Suppose by contradiction  that there are
$u \in D(x) \setminus D(y)$ and $v \in D(y) \setminus D(x)$. Since $L$ is a linear extension of $P$, $I$ is an initial segment of $P$. Since  $u, v\in I$, if $v\leq _P u$ then $v \in D(x)$ and if $u\leq _P v$ then
$u \in D(y)$, a contradiction; so $u \parallel _P v$.
With no loss of generality, we may suppose  $u<_{L} v$. Since $u\leq _Pv$, if $v<_Lx$ then, since $L$ is non-separating, we have $v\leq_P x$, a contradiction. Hence $x<_L v$.  Since $I$ is an initial segment of $L$ we have $x\in L$, contradicting the hypothesis that $x\not \in I$. 
 
(2)  Let $v\in D(y)\setminus D(x)$. Necessarily, $v <_P y$ and, since $I$ an initial segment of $L$,  $v \parallel _P x$.
If
$x\parallel_{P}y$ and
$x <_L y$ then, since  $L$ is non-separating, we have $x <_L v$, which contradicts the fact that $I$  is an initial segment of
$L$.  Hence, either $ x \not\parallel_{P} y$, in which case $x<_{P}y$ or  $y< _{L} x$. \end{proof}

\subsection{The dyadic tree}In the two-dimensional case  we have:

\begin{proposition} If $T_2$ is embeddable in a product of two chains
then both are non scattered. \end{proposition}
We deduce this  from the following proposition. 
\begin{proposition}\label{lem:separdyadic} Every non-separating extension of the dichotomic
tree $T_2$ has order type $\omega (1+ \eta)$.
\end{proposition} 
\begin{proof} We use the \emph{condensation method} (see \cite{rosenstein} pp. 71).
Let $\mathcal L$ be a non-separating linear extension of the order on $T_2$ and $L$ be the corresponding chain. Two elements $x,y\in T_2$ are \emph{equivalent} if the interval they determine in $L$ is finite. This is an equivalence relation. Each classe being an interval of $L$. The set of these equivalences is naturally ordered and the chain $L$ is the lexicographical sum of these equivalence classes. 

\begin{claim}\label{claim:equiv1}Each equivalence class is a subchain of $T_2$ and has order type $\omega$. \end{claim}
\noindent{\bf  Proof of Claim \ref{claim:equiv1}.} We observe that for every $x\in T_2$, one of the two covers  of $x$ in $T_2$, namely $x0$ and $x1$,  is a cover of $x$ in $L$. Indeed, suppose for an example $x0<_Lx1$. If $x0$ is not a cover of $x$ in $L$ there is some $y$ with $x<_L y<_Lx0$. With respect to $T_2$ this element  $y$ is incomparable to  $x$ and  $x0$ (if $y$ was comparable  to  $x$ we would have $x_1\leq_{T_2} y$, whereas if  $y$ was incomparable to 
$x0$, then since $T_2$ is a tree, $y$ would be comparable to $x$). Since $x<_{T_2} x0$, $\mathcal L$ is a separating linear extension, contradicting our hypothese. From this observation and the fact that $\downarrow x$ is finite for every $x\in T_2$, the claim follows.  
\endproof
\noindent\begin{claim} \label{claim:equiv2}The set $D$ of equivalence classes has order type $1+\eta$.
\end{claim}
\noindent{\bf  Proof of Claim \ref{claim:equiv2}.} Since $T_2$ has a least element, $D$ too. Also $D$ has no largest element. Otherwise, let $X$ be the largest class. Pick $x\in X$. Since $X$ is a subchain  of $T_2$, one of the two covers of $x$ is not in $X$; its equivalence class is larger than $X$, a contradiction. Finally, not class $X$ has a cover in $D$. Otherwise, if $Y$ is a cover of $X$, let $y$ be the least element of $Y$. Since $\downarrow y$ is finite, there is $x\in X$ wich is incomparable to $y$ (w.r.t. $T_2$). Let $x'$ be the cover of $x$ in $T_2$ which does not belong to $X$. We have $x<_{T_2}x'$,  $x<_{L}y<_{L}<x'$,  $y$ incomparable to $x$ and $x'$ (w.r.t. $T_2$). This  contradicts the fact that $\mathcal L$ is a non-separating extension. \endproof

 With these claims, the proof of Proposition \ref{lem:separdyadic} is complete.
\end{proof}

\subsection{Non-separating scattered extensions }

The  "bracket relation" 
\begin{equation}\label{fred}
\eta \rightarrow [\eta 
]_{2}^{2}
\end{equation}
a famous unpublished result of F.Galvin,  asserts that if the pairs  of rational numbers are divided into finitely many classes  then there
is a  subset $X$ of the rationals  
which is order-isomorphic to the rationals and such 
that all pairs being to the union of two classes (for a proof, see \cite{todorcevic95} Theorem 6.3 p.44 or \cite {fraissetr} A.5.4 p.412 and, for a
far reaching generalization, see
\cite {devlin}). 
This  result  expresses in a very 
economical way  what  the partitions of pairs look 
like. Indeed, what it really says is this: 
\begin{theorem}\label{G}Let $[\Q]^2$ be the set of pairs of rational numbers and   $A_{1},
 \dots, A_{n}$ be a partition of $[\Q]^2$. For every order $\leq_{\omega}$ on $\Q$ with order type $\omega$  there is a subset $X$ of $\Q$ of order type 
$\eta$ and  indices $i$ and
$j$  (with possibly $i=j$)  such that all pairs of $X$ on which the 
natural order on $\Q$ and the order $\leq_{\omega}$ coincide belong to $A_{i}$ and  all pairs of $X$ on which 
these two orders disagree belong to $A_{j}$.
\end{theorem}
The proof of Theorem \ref{G} from (\ref{fred}) is immediate: intersect the partition 
$A_{1}, \dots ,A_{n}$ 
with the partition $U$,$V$ associated with the two orders ($U$ being made
of  pairs on which the two orders coincide, and $V$ being made of the 
other pairs) and apply iteratively the bracket relation to the 
resulting partition in order to 
find $X$ whose pairs belong to the unions of two classes.

Partitions, or orders, associated to two linear orderings on the same set, 
like the natural order on the rational numbers and  an  order  of type $\omega$
are called {\it sierpinskizations}. Clearly, $\Omega(\eta)$ is a sierpinskization of $\omega \eta$ with $\omega$, whereas $\Omega(\eta)^*$ is a sierpinskization of $\omega \eta$ and $\omega^*$. These two posets are the  basic sierpinskizations of a non scattered chain with $\omega$ and $\omega^*$. Indeed,  if $\alpha$ and $\alpha'$ are two non scattered countable chains then their  sierpinskization with $\omega$ are equimorphic (see \cite{pouzet-zaguia} Corollary  3.4.2).

From Theorem \ref{G},  we have easily:
\begin{proposition}\label{nonsepar,omega(eta)}
Let $P:=(E, \leq)$ be a poset. If neither  $\eta$, $\Omega(\eta)$ nor $\Omega(\eta)^*$ is embeddable in $P$ then for every  non scattered linear extension  $\mathcal L$ of the order on $P$, and every subset $A\subseteq E$ such that $(A, \mathcal L_{\restriction A})$ has type $\eta$   there is an antichain $A'$ of $P$ which is included in $ A$ and such that  $(A', \mathcal L_{\restriction A'})$ has type $\eta$. 
\end{proposition}
\begin{proof}   Let $A\subseteq E$ such that $(A, \mathcal L_{\restriction A})$ has type $\eta$. Let $A_1$, resp. $A_2$,  be the set of pairs $\{x, y\}$ of $[A]^2$ such that $x$ and $y$ are comparable, resp.  incomparable (w.r.t. the order on $P$). Fix an order $\leq_{\omega}$ of type $\omega$ on $A$. Theorem \ref{G} yields a subset $A'$ of $A$ and $i,j\in \{1,2\}$ such that all pairs of $A$ on which the 
order $\mathcal L$ and the order $\leq_{\omega}$ coincide belong to $A_{i}$ and  all pairs of $A$ on which 
these two orders disagree belong to $A_{j}$. As it is easy to check, the  three cases $i=j=1$, $i=1, j=2$ and  $i=2, j=1$ yield respectively that $P_{\restriction A'}$ is  a chain of  type $\eta$, contains a copy of $\Omega(\eta)$ and  contains a copy of $\Omega(\eta)^*$. Thus these    cases are impossible. The only remaining case $a=j=2$ yields the desired conclusion. 
\end{proof}

 \begin{theorem}   \label{equiv} The following properties are equivalent:
 \begin{enumerate}[{(i)}]
 \item  $P$ is the intersection of two scattered chains.
\item \begin{enumerate}

\item $P$ has a non separative scattered extension 
              and     
 \item Neither $\Omega( \eta 
)$ nor $\Omega^*( \eta)$ is embeddable in $P$. 
\end{enumerate}
\end{enumerate}
\end{theorem}  
\begin{proof} $(i)\Rightarrow (ii)$. Item (ii) (a) follows from Theorem \ref{thm:dim2}. Item (ii)(b) follows from the fact that $\Omega(\eta)$ is an obstruction.  
 $(ii)\Rightarrow (i)$. Let $\mathcal C$ be a   non-separative scattered extension of $P$. Let   $\mathcal C'$  be the complement of $\mathcal C$. To conclude, it suffices to prove that $\mathcal C'$ is scattered.  Suppose that it is not. Apply Proposition  \ref{nonsepar,omega(eta)}  to $P$ and $\mathcal L:= \mathcal C'$. Clearly, neither $\eta$, 
 $\Omega( \eta )$ nor $\Omega^*( \eta )$
is embeddable in $P$. Thus,  there is an antichain $A'$ of $P$  such that $C'_{\restriction A'}$ has type $\eta$. But, since $A$ is an antichain of $P$ and the order on $P$ is the intersection of $\mathcal C$ and $\mathcal C'$, it turns out that $\mathcal C_{\restriction A}$ is the dual of $\mathcal C'_{\restriction A}$, thus $(A', \mathcal C_{\restriction A'})$ has type
$\eta$. This contradicts the fact that $\mathcal C$ is scattered.\end{proof}\\

\noindent{\bf Bibliographical comments.} The posets $T_2$,  $\Omega (\eta)$, $B(\check \eta)$ have been considered previously. 
Pouzet and Zaguia \cite{pouzet-zaguia}
proved  that \emph{the set $\mathcal J(P)$ of ideals of a poset $P$ contains no chain isomorphic to $\eta$
if and only if $P$ contains no chain isomorphic to $\eta$ and
 no subset isomorphic to $\Omega (\eta )$}. In  \cite{duffus-pouzet-rival} it is shown that \emph{if  a poset $P$ contains $B(\check \eta )$, then $N(P)$
contains a chain isomorphic to $\eta$}. In fact, 
$N(B(\check \eta))$ is isomorphic to the disjoint union $\Q\times 2 \cup I(\Q)$ equipped with the following ordering:
\begin{enumerate}
\item $\emptyset \leq (x,0)\leq I\leq (y,1) \leq \Q$ for  $x\in I\subseteq (\leftarrow y[$, with $I\in I(\Q)$, $y\in\Q$.
\item $I\leq J$ if $I\subseteq J$ and $I,J\in I(\Q)$.
\end{enumerate}
In \cite{pousikadzag} it is shown that  \emph{the class  of posets whose MacNeille completion is scattered is characterized by eleven obstructions}. One can check that obstructions distinct from  $\eta$ and  $B(\check \eta)$ do no yield interesting obstructions to $\mathcal L_{\mathcal S}(<n)$. 
In an unpublished paper with E.C.Milner \cite{milner}  it is shown that if {the set $\mathcal J(P)$ of ideals of a poset $P$ is topologically closed in $\mathfrak P(P)$, it is topologically scattered if and only if it is order scattered and the binary tree $T_2$ is not embeddable in $P$}. From this follows that \emph{an algebraic lattice $T$ is topologically scattered if and only if it is order scattered and neither $T_2$ nor $\Omega (\eta )$ are embeddable in the join-semilattice of compact elements of $T$}. In contrast, we may note that \emph{an algebraic distributive lattice is topologically scattered if and only if it is order scattered}, an important result due to Mislove \cite{mislove}.

\end{document}